\documentclass[12pt]{amsart}
\usepackage{mathtools, amssymb, bm, amsthm, amsmath, mathrsfs}

\usepackage{geometry}
\usepackage[utf8]{inputenc}

\usepackage{scalerel}[2016-12-29]
\def\stretchint#1{\vcenter{\hbox{\stretchto[440]{\displaystyle\int}{#1}}}}
\def\bs{\mkern-12mu} 

\usepackage{bigints}

\usepackage{biblatex}
\usepackage{color}

\usepackage{hyperref}

\newtheorem{theorem}{Theorem}
\newtheorem{proposition}[theorem]{Proposition}
\newtheorem{lemma}{Lemma}[section]

\theoremstyle{definition}
\newtheorem{remark}{Remark}[section]

\numberwithin{theorem}{section}
\numberwithin{equation}{section}

\newcommand{\R}{\mathbb R}
\newcommand{\Z}{\mathbb Z}
\newcommand{\qZ}{(\Z/q\Z)^\times}
\newcommand{\N}{\mathbb N}
\newcommand{\Q}{\mathbb Q}
\newcommand{\C}{\mathbb C}

\newcommand{\eps}{\varepsilon}

\begin{document}

\title[variance of square-full integers]{variance of square-full integers in short intervals and arithmetic progressions}

\author{Yotsanan Meemark}
\author{Watcharakiete Wongcharoenbhorn}

\address{Department of Mathematics and Computer Science, Faculty of Science, Chulalongkorn University, Bangkok, Thailand 10330}

\email{yotsanan.m@chula.ac.th}
\email{w.wongcharoenbhorn@gmail.com}


\begin{abstract}
Define a natural number $n$ as a \textit{square-full} integer if for every prime $p$ such that $p|n$, we have $p^2|n$. In this paper, we establish an upper bound on the variance of square-full integers in short intervals of an expected order, under the assumption of a certain quasi-Riemann hypothesis. We also prove an asymptotic formula for the variance in arithmetic progressions, averaging over a quadratic residue and a nonresidue by a half, which is of smaller order of magnitude than the aforementioned bound for all primes $q\gg x^{51/114+\eps}$.
\end{abstract}

\maketitle

\begin{section}{Introduction}
A positive integer $n$ is called \textit{square-full} if for every prime divisor $p$ of $n$, we have $p^2\mid n$. As conventional, we denote $\mu(n)$ by the M\"{o}bius function on $n$. An important property of square-full integers is that they can be written as $a^2b^3$ uniquely for $b$ square-free. Let $Q(x)$ be the number of square-full integers not exceeding $x$. Bateman and Grosswald showed in \cite{Bateman-Grosswald} that there exists a constant $c>0$ such that
\begin{align}
\label{eq:counting-sqfull}
Q(x):&=\sum_{a^2b^3\leqslant x}\mu^2(b)=\dfrac{\zeta(\frac{3}{2})}{\zeta(3)}x^{\frac{1}{2}}+\dfrac{\zeta(\frac{2}{3})}{\zeta(2)}x^{\frac{1}{3}}+O(x^{\frac{1}{6}}\exp(-c\log^{4/7}x(\log\log x)^{-3/7})),
\end{align}
which improved upon the result of Erd\H{o}s and Szekeres \cite{square-full-Erdos-Szekeres} by incorporating the term $x^{\frac{1}{3}}$ explicitly with the presented error term. Counting square-full integers is very closely related to the \textit{$(2,3)$-divisor problem}, which determines for all $x\geqslant 1$, the smallest possible order of the error term $\Delta_{2,3}(x)$ occurring from the formula
\[
\sum_{a^2b^3\leqslant x}1=\zeta(\tfrac{3}{2})x^{\frac{1}{2}}+\zeta(\tfrac{2}{3})x^{\frac{1}{3}}+\Delta_{2,3}(x).
\]

It is possible to show that $\Delta_{2,3}(x)=O(x^{\frac{2}{15}})$ via the Dirichlet hyperbola method using the bounds on exponential sums. Let $\eps>0$ be a small positive number. Obtaining an upper bound of $\Delta_{2,3}(x)$ of order as small as that expected in the work of Kr\"{a}tzel \cite{Kratzel-theorem} of order $O_\eps(x^{\frac{1}{10}+\eps})$ seems to be very difficult. In fact, he proved that $\Delta_{2,3}(x)=\Omega(x^{\frac{1}{10}})$, which means that for some constant $c>0$, the absolute value of this term can be larger than $cx^{\frac{1}{10}}$ infinitely often. In addition, he proved similar results for the $(\alpha,\beta)$-divisor problem, defined similarly on positive integers $\alpha,\beta$ including the original Dirichlet divisor problem. Subsequently, Balasubramanian, Ramachandra, and Subbarao \cite{error-term-in-square-full} established an analogous result to Kr\"atzel on square-full integers. In particular, they proved a similar Big Omega result that $\Delta(x)=\Omega(x^{\frac{1}{10}})$, where $\Delta(x)$ denotes the error term in Equation \eqref{eq:counting-sqfull}. It seems reasonable to conjecture that the error term in counting square-full numbers is also of the order $O_\eps(x^{\frac{1}{10}+\eps})$ for all $\eps>0$. In \cite{Bateman-Grosswald}, Bateman and Grosswald also indicated that the result, such as $\Delta(x)=O_\eps(x^{\frac{1}{6}-\delta+\eps})$ for any small $\delta>0$, would inevitably prove a certain quasi-Riemann hypothesis (QRH). Hence, obtaining this expected order is certainly difficult. The known results on the distribution of square-full numbers of small error terms are dependent on the Riemann hypothesis (RH). In fact, with current methods, it seems that the order of $x^{1/8}$ is the limit of the full strength of the RH. The best result in this direction is due to Wang \cite{Wang} showing that $\Delta(x)=O(x^{328/2333})$ under the RH. Here, $328/2333\approx 0.1406$.\\

Not only is the literature concerned with counting square-full integers in long intervals, such as $[1,X]$, they are also interested in counting them in short intervals $(x,x+H]$ for $x\in[X,2X]$. The average along such an interval is deemed short since $H$ is much smaller than $X$ in the sense that $H=o(X)$. For later references, we denote $1^{\blacksquare}(n)$ by the characteristic function of square-full integers. Consequently, the Erd\H{o}s-Szekeres theorem, i.e., \eqref{eq:counting-sqfull} implies that
\begin{equation}
\label{eq:sqfull-inshort-interval}
\sum_{x<n\leqslant x+H}1^\blacksquare(n)\sim\dfrac{\zeta(\frac{3}{2})H}{2\zeta(3)\sqrt{x}},
\end{equation}
provided $H\geqslant x^{\frac{1}{2}+\frac{1}{6}}$. The number $\frac{1}{6}$ is reduced subsequently to extend the possible range that \eqref{eq:sqfull-inshort-interval} still holds. Given that the error term is conjectured to be of the order $O_\eps(x^{\frac{1}{10}+\eps})$, this number could be replaced by $\frac{1}{10}$ or even lower. The best-known result in this direction is due to Trifonov \cite{Trifonov}, who showed that we can take $\frac{19}{154}\approx 0.1233$ in place of $\frac{1}{6}$. 

Determining the error term appearing in the asymptotic formula of $\sum_{x<n\leqslant x+H} 1^\blacksquare(n)$ is more difficult than the previous long average and is almost certainly deeper than the $\frac{1}{10}$-QRH. Here, we mean $r$-QRH for some $r\in[0,\frac{1}{2})$ by the conjecture that asserts that the zeros $\rho$ of the zeta function satisfy $\Re\rho\in[\frac{1}{2}-r,\frac{1}{2}+r]$. Whence, $0$-QRH corresponds to the full RH. In harmony with this, it is reasonable to consider working on the variance of the error term as initial steps towards the more challenging problems. Let us write $\mathscr D_r(x;y):=(x+y)^r-x^r$ for any positive $r\in\Q$, which serves as a difference between the main terms that square-full numbers are counted in short intervals. Chan \cite{Chan1} and \cite{Chan2} studied the expression
\[
\stretchint{7ex}_{\bs X}^{2X}\left|\sum_{x<n\leqslant x+y} 1^\blacksquare(n)-\dfrac{\zeta(\frac{3}{2})}{\zeta(3)} \mathscr D_{\frac{1}{2}}(x;y)\right|^2dx
\]
in which $y:=2\sqrt{x}H+H^2$. He proved nice results in the papers, including its asymptotic formula. Particularly, he unconditionally proved for $X^\eps\leqslant H\leqslant X^{\frac{189}{1046}-\eps}$ the asymptotic form of the above variance. Under the Lindel\"{o}f hypothesis, he could also extend the range to $H\leqslant X^{\frac{1}{4}-\eps}$, which is the best possible if the secondary main term is not included in the calculation.

\subsection{Main results.}
Our work is inspired by the breakthrough work of Gorodetsky, Matom\"{a}ki, Radziwi\l\l\text{ }and Rodgers \cite{GoroMato} studying the variance in square-free integers in short intervals and arithmetic progressions, as well as the aforementioned pioneering works \cite{Chan1} and \cite{Chan2} of Chan on square-full integers in short intervals. In particular, we modify slightly the method in \cite{GoroMato} to prove our two main theorems, Theorem \ref{thm:QRH-equiv-variance-sqfull} and Theorem \ref{thm:sqfull-in-arith} below.

\begin{theorem}
\label{thm:QRH-equiv-variance-sqfull}
The $\frac{1}{10}$\textup{-QRH} is equivalent to the following statement. For $\eps\in(0,\frac{1}{1000})$, we have for $X^{\eps}\leqslant H\leqslant X^{1-\eps}$ that
\begin{equation}
\label{eq:Variance-of-full-short-ranges}
\dfrac{1}{X}\stretchint{7ex}_{\bs X}^{2X}\left|\sum_{\substack{x<n\leqslant x+H \\ }}1^{\blacksquare}(n)-\left(\dfrac{\zeta(\frac{3}{2})}{\zeta(3)} \mathscr D_{\frac{1}{2}}(x;H)+\dfrac{\zeta(\frac{2}{3})}{\zeta(2)}\mathscr D_{\frac{1}{3}}(x;H)\right)\right|^2dx\ll_{\eps} X^{\frac{1}{5}+5\eps}.
\end{equation}
\end{theorem}

With the steps used in proving Theorem \ref{thm:QRH-equiv-variance-sqfull}, we can also prove the analogous result for $(\alpha,\beta)$-full integers, defined naturally from the $(\alpha,\beta)$-divisor problem, for $\alpha,\beta>1$. Indeed, it is attributed in \cite{Suryanarayana-generalized-sqfull} by Suryanarayana that Cohen introduced such numbers as early as 1963. He also proved the asymptotic relation of them and called them generalized square-full integers. Thus, the RH is equivalent to all of these variance results of such $(\alpha, \beta)$-full integers combined. It also seems to the authors that the method that could allow us to prove the variance for square-full integers should be adapted to any $(\alpha,\beta)$-full integers. In this regard, the $\frac{1}{10}$-QRH is very close to the full RH, which is pretty strange.

As mentioned in \cite{GoroMato}, the variance of this type is difficult for large $H=o(X)$ since there are a number of interferences between each integrand. We see that the result over merely short ranges of large $H$, e.g., for $H\in [X^{1-2\eps}, X^{1-\eps}]$, is sufficient to imply the $\frac{1}{10}$-QRH.\\

It is quite reasonable to expect that, for $0\not=|h|\leqslant x^{1-\eps}$,
\begin{align*}
    \sum_{n\leqslant x} 1^{\blacksquare}(n)1^{\blacksquare}(n+h)\sim C_h\log x,
\end{align*}
for some constant $C_h>0$ depending only on $h$. This is stronger than the consequence of the \textit{abc} conjecture, which the authors are able to prove. One can show that the left-hand side above is $\ll _{\eps,h} x^{\eps}$ assuming the $abc$ conjecture. If we assume further, from this weakened result, that the implied constant is uniform on $|h|\leqslant H$, which we are unable to show, it is possible to recover \eqref{eq:Variance-of-full-short-ranges} but only in the range $H\leqslant X^{\frac{27}{40}}$ using the interchange of the sums and integrals.
\\

\smallskip

There has also been interest in counting square-full integers in arithmetic progressions. Fix $q\geqslant 1$ be an integer and $\ell\in(\Z/q\Z)^\times$. Sometimes, we refer to this problem as a $q$-aspect in counting square-full numbers. We define $Q(x;q,\ell)$ as a number of square-full integers not exceeding $x$ that are congruent to $\ell$ modulo $q$. Similar to Equation \eqref{eq:counting-sqfull} we have the result of Chan \cite{Chan-arith}
\begin{align}
\label{eq:counting-sqfull-in-arith}
    Q(x;q,\ell):=\sum_{\substack{a^2b^3\leqslant x\\ a^2b^3\equiv\ell(q)}}\mu^2(b)=\dfrac{A_{q,\ell}}{\zeta(3)}\dfrac{x^{\frac{1}{2}}}{q}+\dfrac{B_{q,\ell}}{\zeta(2)}\dfrac{x^{\frac{1}{3}}}{q}+O_\eps\left(x^{\frac{1}{6}}q^{\frac{1}{12}+\eps}+\frac{x^{\frac{1}{5}}}{q^{\frac{1}{5}-\eps}}\right),
\end{align}
where $A_{q,\ell},B_{q,\ell}$ are certain constants. This improved the error term from the work of Chan and Tsang in \cite{chan-tsang} and is currently the best-known result. In this paper, we also study the variance in this direction and prove the following theorem.

\begin{theorem}
\label{thm:sqfull-in-arith}
    Let $\eps\in(0,\frac{1}{1000})$ be given, and let $q$ be a prime. Let $\alpha$ be a quadratic nonresidue modulo $q$. We have for all $
    x^{\frac{51}{114}+\eps}\leqslant q\leqslant x^{1-50\eps}$,
    \begin{align}
    \label{eq:main-eq-them-arith}
        \dfrac{1}{\varphi(q)}\sum_{\ell\in (\Z/q\Z)^\times} & \left| \frac{1}{2}\sum_{\substack{n\equiv \ell,\alpha\ell(q)\\ x<n\leqslant 2x}} 1^{\blacksquare}(n)-\left(\dfrac{\zeta(\frac{3}{2})}{\zeta(3)}\frac{1}{q}\left(1-\frac{1}{q}\right)\mathscr D_{\frac{1}{2}}(x;x)+\dfrac{\zeta(\frac{2}{3})}{\zeta(2)}\frac{1}{q}\left(1-\dfrac{1}{q^{\frac{2}{3}}}\right)\mathscr{D}_{\frac{1}{3}}(x;x)\right)
         \right|^2 \nonumber \\
         &\sim C(x/q)^{\frac{1}{6}},
    \end{align}
    where $$C:= \frac{\zeta(\frac{11}{6})}{2}\prod_p\left(1-\dfrac{1}{p^2}+\dfrac{2}{p^{6}}+\dfrac{2}{p^7}\right)\stretchint{7ex}_{\bs 0}^{\infty}|W(y)|^2y^{\frac{5}{6}}dy,$$
    with a smooth function $W:=\widehat{f1_{[1,2]}}$ where $f$ is defined as in Equation \eqref{eq:approx=characteristic-sqfull} and $1_{[1,2]}$ is the characteristic function of $[1,2]$.
\end{theorem}
\begin{remark}
We note here that the main term in the variance is merely an average of the two possible primary main terms of order $x^{1/2}/q$. In Chan \cite{Chan-arith}, the constant $A_{q,\ell}$ in Equation \eqref{eq:counting-sqfull-in-arith} is given by
$$\prod_{p|q}(1-p^{-3})^{-1}\sum_{b\geqslant 1} \dfrac{N_2(\ell b;q)}{b^{3/2}}\dfrac{x^{1/2}}{q},$$
where $N_2(n;q) := \#\{x\in (\Z/q\Z)^{\times}: x^2\equiv n\pmod{q}\}$. This number depends on whether $\ell$ is a quadratic residue modulo $q$ or not. Additionally, our result is slightly better than the trivial result, which occurs when we bound each summand by $x^{1/2}/q$, the left-hand side of Equation \eqref{eq:main-eq-them-arith} is $\ll (x/q)^{1/6}$ when $q\gg x^{5/11}$.

Since the summand is of order $(x/q)^{1/12}$ on average over $q\gg x^{51/114+\eps}$ while the expected order of $\Delta(x)$ is $O_\eps(x^{1/10+\eps})$, we believe that the range of $q$ in Theorem \ref{thm:sqfull-in-arith} might not be possible to extend to the whole range $q\geqslant 1$.
\end{remark}

\subsection{Fractional Brownian motion.}
The statistical results of certain integer sequences, like the one we presented, always (or at least are believed to) deeply connect to some well-known distributions. For instance, the map such as
\begin{align}
\label{proc:primes}
t\mapsto \frac{1}{H^{1/2}} \sum_{x<p\leqslant x+tH} \left(\log p-1\right)
\end{align}
where the index $p$'s are primes, is conjectured to behave like the usual Brownian motion, i.e., that the process will tend to the fractional Brownian motion of Hurst parameter $\frac{1}{2}$. In Section 1.2 of \cite{GoroMato}, they also discussed the features involving square-free integers and fractional Brownian motion, in which it seems to involve the Hurst parameter $\frac{1}{4}$. Therefore, we are particularly interested in interpreting these expected phenomena for square-full integers as well. We think that we should consider the process
\begin{align}
\label{proc:square-fulls}
t\mapsto \frac{1}{H^{3/5}} \sum_{x<n\leqslant x+tH} \left(\frac{1^{\blacksquare}(n)}{\frac{\zeta(\frac{3}{2})}{2\zeta(3)n^{\frac{1}{2}}}+\frac{\zeta(\frac{2}{3})}{3\zeta(2)n^{\frac{2}{3}}}}-1 \right)
\end{align}
which appears to correspond to the Hurst parameter $\frac{3}{5}$. We illustrate the evaluations of the map \eqref{proc:primes} and \eqref{proc:square-fulls} below by setting $x=5\times 10^9, H=46,674,434$ and $t$ ranges between $0$ and $1$ in both figures.\\

\begin{figure}[ht]
    \centering
    \includegraphics[width=0.8\linewidth]{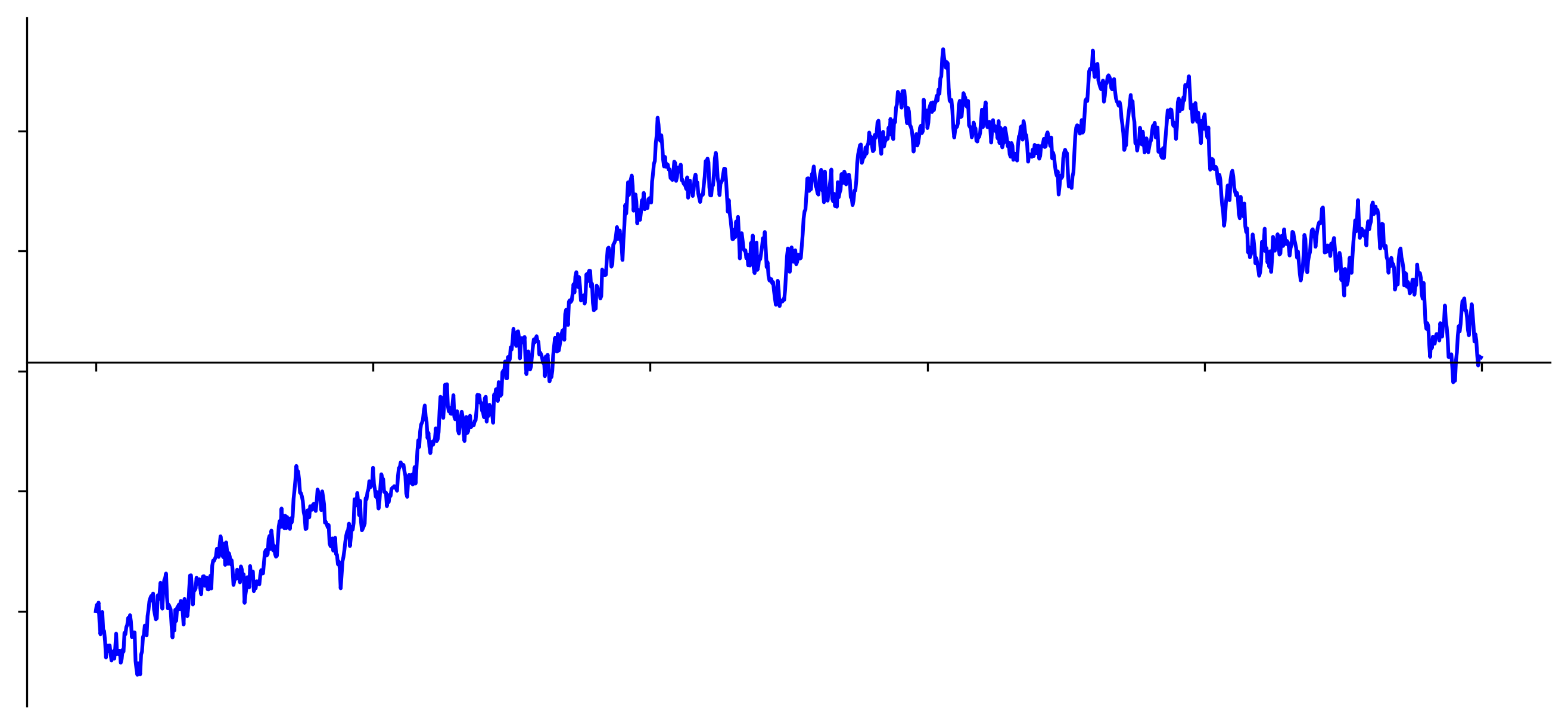}
    \caption{Prime numbers: Values of \eqref{proc:primes} for $x=5\times 10^9, H=46,674,434$ and $t$ ranges between $0$ and $1$.}
    \label{fig:Prime numbers}
\end{figure}

\begin{figure}[ht]
    \centering
    \includegraphics[width=0.8\linewidth]{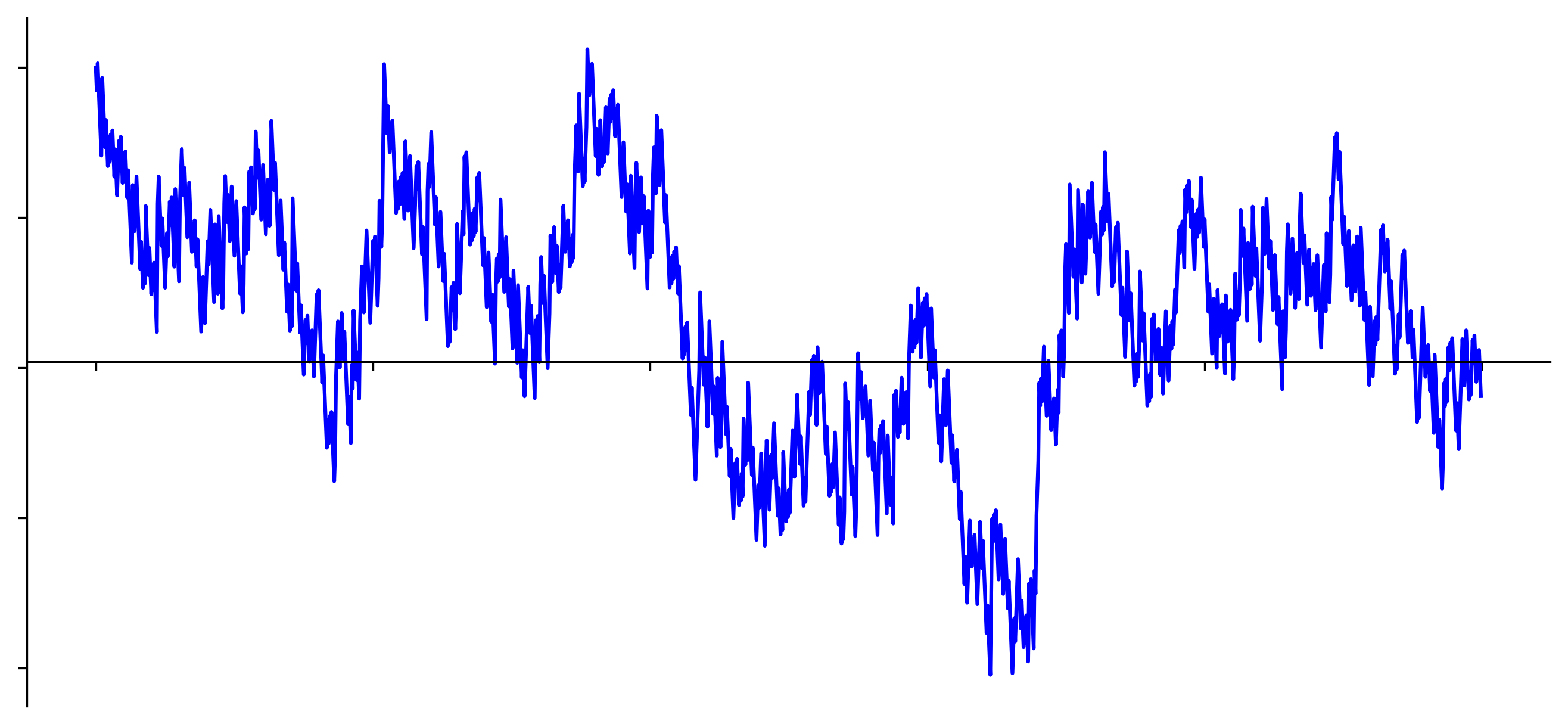}
    \caption{Square-full integers: Values of \eqref{proc:square-fulls} for $x=5\times 10^9, H=46,674,434$ and $t$ ranges between $0$ and $1$.}
    \label{fig:Square-full integers}
\end{figure}

Our paper is organized as follows. Firstly, the variance of the $(2,3)$-divisor problem and of square-full integers in short intervals will be discussed in Section \ref{sec:Variance of $(2,3)$-divisor problem and of square-full integers in short intervals}, where we establish Theorem \ref{thm:QRH-equiv-variance-sqfull}. Then, Section \ref{sec:Propositions in the $q$-aspect and Variance of square-full integers in arithmetic progressions} alone is dedicated to conferring about propositions related to the variance in arithmetic progressions, i.e., the variance in the $q$-aspect. This section culminates in the proof of Theorem \ref{thm:sqfull-in-arith}.
\end{section}
\section{Variance of \texorpdfstring{$(2,3)$}{(2,3)}-divisor problem and of square-full integers in short intervals}
\label{sec:Variance of $(2,3)$-divisor problem and of square-full integers in short intervals}

This section is devoted to the proof of Theorem \ref{thm:QRH-equiv-variance-sqfull}. Firstly, we deduce that the result on the variance of square-full integers in short intervals of expected order implies the $\frac{1}{10}$-QRH. Then, we state and prove the result involving the variance of the $(2,3)$-divisor problem. This will be employed in proving, under the $\frac{1}{10}$-QRH, that the variance of square-full integers in short intervals is indeed of the expected order.

\subsection{Result on the variance implies $\frac{1}{10}$-QRH}

Suppose that the $\frac{1}{10}$-QRH is false, so there is a zero $\rho$ of the zeta function such that $\Re\rho>\frac{3}{5}$. Let $6\Theta>\frac{3}{5}$ be the supremum of the real part of zeroes of $\zeta$ and let $\delta>0$ be small such that $\Theta-\frac{\delta}{2}\geqslant \frac{1}{10}+20\delta$.\\

Since any square-full integer can be written uniquely as $a^2b^3$ with $b$ square-free, we can count square-full numbers $n=a^2b^3$ using $\mu^2(b)$. Thus, we have from the Cauchy-Schwarz inequality with \eqref{eq:Variance-of-full-short-ranges} that
\begin{align*}
    \frac{1}{X}&\stretchint{7ex}_{\bs X}^{2X} \Phi\left(\frac{x}{X}\right)\left(\sum_{x<a^2b^3c^6\leqslant x+H} \mu(c)-\frac{\zeta(\tfrac{3}{2})}{\zeta(3)}\mathscr{D}_{\frac{1}{2}}(x;H)-\frac{\zeta(\tfrac{2}{3})}{\zeta(2)}\mathscr{D}_{\frac{1}{3}}(x;H)\right) dx
    \\ &\ll_\Phi \left(\frac{1}{X}\stretchint{7ex}_{\bs X}^{2X}\left| \sum_{x<n\leqslant x+H} 1^\blacksquare(n)-\frac{\zeta(\tfrac{3}{2})}{\zeta(3)}\mathscr{D}_{\frac{1}{2}}(x;H)-\frac{\zeta(\tfrac{2}{3})}{\zeta(2)}\mathscr{D}_{\frac{1}{3}}(x;H)\right|^2dx\right)^\frac{1}{2}\ll_\delta H^{\frac{1}{10}+10\delta}
\end{align*}
for any compactly supported smooth function $\Phi(x)$ with $\text{supp}(\Phi)\subset [1,2]$. Here, we have used the identity $\mu^2(b)=\sum_{d|b}\mu(d)$. From the above calculation, we obtain
\begin{align}
\label{eq:imply(1/10)-QRH after CS}
    \dfrac{1}{X}\stretchint{7ex}_{\bs X}^{2X} \Phi\left(\frac{x}{X}\right) \sum_{x<a^2b^3c^6\leqslant x+H} \mu(c)=\frac{\zeta(\tfrac{3}{2})}{2\zeta(3)}X^{\frac{1}{2}}\Psi_{\tfrac{H}{X}}(\tfrac{1}{2})+\dfrac{\zeta(\tfrac{2}{3})}{3\zeta(2)}X^{\frac{1}{3}}\Psi_{\tfrac{H}{X}}(\tfrac{1}{3})+O_\delta(H^{\frac{1}{10}+10\delta}),
\end{align}
where we define for any complex number $s$,
\begin{align*}
    \Psi_{\tfrac{H}{X}}(s) := \frac{1}{s}\stretchint{7ex}_{\bs\R} \left(\left(x+\frac{H}{X}\right)^s-x^s\right)\Phi(x) dx,
\end{align*}
and we shift the term of the form $x+H$ in the integral to get this definition. Now, we employ Perron's formula on the left-hand side sum of \eqref{eq:imply(1/10)-QRH after CS} and interchange the integrals to obtain that the left of \eqref{eq:imply(1/10)-QRH after CS} is
\begin{align*}
    \dfrac{1}{2\pi i}\stretchint{7ex}_{\bs(2)} \dfrac{\zeta(2s)\zeta(3s)}{\zeta(6s)} \cdot X^s\Psi_{\tfrac{H}{X}}(s) ds.
\end{align*}
Let $\widetilde{f}(s) := \int_{\R^+} f(x)x^{s-1} dx$ be the Mellin transform of $f$. As in the proof of Theorem 3 in \cite{GoroMato}, our first step is to show that
\begin{align}
\label{eq:(1/10)implication-firststep}
    \stretchint{7ex}_{\bs (\sigma)} \dfrac{\zeta(2s)\zeta(3s)}{\zeta(6s)}\cdot X^s\widetilde{\Phi}(s) ds \ll_\delta X^{\Theta-\frac{\delta}{3}},
\end{align}
where we move the contour from the path $(2)$ to $(\sigma)$ for which $\sigma=\Theta+\frac{\delta}{10}$ to avoid any zeroes of $\zeta(6s)$. The poles from $s=\frac{1}{2},\frac{1}{3}$ correspond to the main term on the right-hand side of \eqref{eq:imply(1/10)-QRH after CS}. Thus, we are left to replace $\Psi_{\tfrac{H}{X}}(s)$ by $\widetilde{\Phi}(s)$.\\

Let us first proof \eqref{eq:(1/10)implication-firststep}. For any compactly supported smooth function $\vartheta$ such that $\text{supp}(\vartheta)\subset [\frac{1}{1000},1000]$, and $A\in\N$, we see that 
\begin{align*}
    \stretchint{7ex}_{\bs \R}\vartheta(x) x^sdx=\frac{(-1)^A}{(s+1)(s+2)\ldots (s+A)} \stretchint{7ex}_{\bs\R} \vartheta^{(A)}(x)x^{s+A}dx\ll_{A} (1+|\Im(s)|)^{-A},
\end{align*}
by repeatedly using Taylor's theorem on $x^sdx$. Then, we have by Taylor's expansion and the above that
\begin{align*}
    &\Psi_{\tfrac{H}{X}}(s) =\frac{1}{s}\stretchint{7ex}_{\bs\R} x^s\left(\Phi'(x)\left(-\frac{H}{X}\right)+\sum_{j\geqslant 2} \dfrac{\Phi^{(j)}(x)}{j!}\left(-\frac{H}{X}\right)^j\right) dx
    \\ &= \left(-\frac{H}{X}\right)\stretchint{7ex}_{\bs \R} \frac{x^s}{s}d\Phi(x)+O_A\left(\sum_{j\geqslant 2} \frac{X^{-\delta}(1+\Im(s))^{-A}}{j!}\right)=\left(\frac{H}{X}\right)\widetilde{\Phi}(s)+O_A\left(\frac{X^{-\delta}}{(1+|\Im(s)|)^{A}}\right).
\end{align*}
We will combine with the estimate
\begin{align*}
    \dfrac{\zeta(2s)\zeta(3s)}{\zeta(6s)}\ll_\delta (1+|\Im(s)|)^{C}
\end{align*}
for some constant $C=C(\delta)$, and the inequality holds uniformly for all $\Re s\geqslant \sigma$. Hence,
\begin{align*}
    \stretchint{7ex}_{\bs (\sigma)}\dfrac{\zeta(2s)\zeta(3s)}{\zeta(6s)}\cdot X^s\Psi_{\tfrac{H}{X}}(s) ds = \frac{H}{X}\stretchint{7ex}_{\bs (\sigma)} \dfrac{\zeta(2s)\zeta(3s)}{\zeta(6s)}\cdot X^s\widetilde{\Phi}(s) ds+O_\delta\left(X^{\sigma-\delta}(1+\Im(s))^{C-A}\right),
\end{align*}
where the last error term is $\ll_\delta X^{\Theta-\frac{\delta}{2}}$ for any sufficiently large $A\geqslant C$, by the definition of $\sigma$. Equation\eqref{eq:(1/10)implication-firststep} now holds since the remaining error term $(X/H)O_\delta(H^{\frac{1}{10}+10\delta})\ll H^{\frac{1}{10}+11\delta}\leqslant X^{\Theta-\frac{\delta}{2}}$ by the definition of $\delta$.
\\

We shall now prove that
\begin{align}
\label{eq:Mellin's discrete convolution}
    \sum_{n\geqslant 1} a(n):=\sum_{n\geqslant 1} 1^{\blacksquare}(n)\Phi(n/x)=\frac{1}{2\pi i}\stretchint{7ex}_{\bs (2)} \dfrac{\zeta(2s)\zeta(3s)}{\zeta(6s)}\cdot x^s\widetilde{\Phi}(s) ds.
\end{align}
Typically, this is Mellin's discrete convolution; however, in our situation, this is equivalent to Perron's formula. We have that the partial summation (up to $y$) of the left of \eqref{eq:Mellin's discrete convolution} is, by partial summation,
\begin{align*}
    \left(\sum_{n\leqslant y} 1^\blacksquare(n)\right)\Phi(y/x)-\stretchint{7ex}_{\bs 1}^y \left(\sum_{n\leqslant u} 1^\blacksquare(n)\right)\dfrac{\partial\Phi(u/x)}{\partial u} du
\end{align*}
and the first term vanishes as $y\to\infty$. By using Perron's formula in the above and then partial integration, the left-hand side of \eqref{eq:Mellin's discrete convolution} is
\begin{align*}
    -\frac{1}{2\pi i} \stretchint{7ex}_{\bs (2)}\dfrac{\zeta(2s)\zeta(3s)}{\zeta(6s)}\stretchint{7ex}_{\bs 1}^\infty \frac{u^s}{s}\dfrac{\partial\Phi(u/x)}{\partial u}du ds =\frac{1}{2\pi i}\stretchint{7ex}_{\bs (2)} \dfrac{\zeta(2s)\zeta(3s)}{\zeta(6s)}\stretchint{7ex}_{\bs \R^+} \Phi(u/x)u^{s-1}du,
\end{align*}
which is the right-hand side of \eqref{eq:Mellin's discrete convolution} by changing variable $u\mapsto ux$.\\ 

Let the left-hand side of Equation \eqref{eq:(1/10)implication-firststep} be $A(X)$, we consider, for a fixed $s\in\C$,
\begin{align*}
    \stretchint{7ex}_{\bs 0}^\infty \frac{A(x)}{x^{s+1}} dx &= \stretchint{7ex}_{\bs 0}^{\infty} \dfrac{1}{x^{s+1}}\left(\sum_{n\geqslant 1} a(n)-\dfrac{\zeta(\tfrac{3}{2})\widetilde{\Phi}(\tfrac{1}{2})}{2\zeta(3)}x^{\frac{1}{2}}-\dfrac{\zeta(\tfrac{2}{3})\widetilde{\Phi}(\tfrac{1}{3})}{3\zeta(2)}x^{\frac{1}{3}}\right) dx
\end{align*}
since we evaluate the poles at $s=\frac{1}{2}, \frac{1}{3}$ in the integral of \eqref{eq:Mellin's discrete convolution}. Let $\Re s>1$, by the change of variable $x\mapsto nu^{-1}$ the first term on the right is
\begin{align*}
    \sum_{n\geqslant 1} 1^\blacksquare(n)\stretchint{7ex}_{\bs 0}^\infty x^{-s-1}\Phi(n/x) dx=\sum_{n\geqslant 1}\frac{1^\blacksquare(n)}{n^s}\stretchint{7ex}_{\bs 0}^\infty \Phi(u)u^{s-1} du=\dfrac{\zeta(2s)\zeta(3s)}{\zeta(6s)}\widetilde{\Phi}(s).
\end{align*}
Thus, we obtain for $\Re s>1$,
\begin{align*}
    \stretchint{7ex}_{\bs 1}^\infty \frac{A(x)}{x^{s+1}} dx &= \dfrac{\zeta(2s)\zeta(3s)}{\zeta(6s)}\widetilde{\Phi}(s)+\dfrac{\zeta(\tfrac{3}{2})\widetilde{\Phi}(\tfrac{1}{2})}{2\zeta(3)}\frac{1}{s-\frac{1}{2}}+\dfrac{\zeta(\tfrac{2}{3})\widetilde{\Phi}(\tfrac{1}{3})}{3\zeta(2)}\frac{1}{s-\frac{1}{3}}.
\end{align*}
The left-hand side above can be extended analytically to $\Re s>\Theta-\frac{\delta}{5}$ by \eqref{eq:(1/10)implication-firststep}, but the right-hand side has a pole between $[\Theta-\frac{\delta}{5},\Theta]$, a contradiction. Hence, $\frac{1}{10}$-QRH is true assuming the variance.\\

Now, as promised, we state and prove the unconditional result on the variance of the $(2,3)$-divisor problem. Proposition \ref{prop:var-(2,3)-divisor-problem} below is a special case in the work of Ivi\'{c} \cite{Ivic-divisor-general} but we will provide an alternative proof for the sake of completeness. Our proof uses Lemma \ref{lem:intboundon_exp(2pi i)} and Lemma \ref{lem:reducetoexp1} below with Voron\"{o}i's argument, which was used in attacking the original divisor problem.\\

As a convention we write $\psi(u):=u-[u]-\tfrac{1}{2}$. The Fourier series of $\psi(t)$ is not absolutely convergent, so we may write for a fixed parameter $J>1$, which will be specified later, its approximation $\widetilde{\psi}(x)$,
$$\widetilde{\psi}(x):=\dfrac{J}{2}\stretchint{7ex}_{\bs -1/J}^{1/J} \psi(x+u) du,$$
which is Lipschitz, and hence, its Fourier series converges absolutely. The error occurred from the approximation $h(x):=|\widetilde{\psi}(x)-\psi(x)|$, is also Lipschitz. The arithmetic function $1_{2,3}(n)$ is the characteristic function of the numbers that are of the form $a^2b^3$ for some integers $a,b$.

\begin{proposition}[Variance of $(2,3)$-divisor problem]
\label{prop:var-(2,3)-divisor-problem}
For arbitrarily small $\eps>0$, we have that
\begin{align}
\label{eq:(2,3)-divisor problem-variance}
\dfrac{1}{X}\stretchint{7ex}_{\bs X}^{2X}\left|
\sum_{\substack{n\leqslant x}}1_{2,3}(n) - \left(\zeta(\tfrac{3}{2})x^{\frac{1}{2}}+\zeta(\tfrac{2}{3})x^{\frac{1}{3}}\right)
\right|^2
dx\ll_\eps X^{\frac{1}{5}+\eps}.
\end{align}
\end{proposition}

\begin{remark}
Since square-full integers correspond to the $(2,3)$-divisor problem, it is natural to inquire about the variance of the $(2,3)$-divisor problem in short intervals. In particular, the appropriate result involving
\begin{align*}
\dfrac{1}{X}\stretchint{7ex}_{\bs X}^{2X} \left|\sum_{x<n\leqslant x+H}1_{2,3}(n)-\left(\zeta(\tfrac{3}{2})\mathscr D_{\frac{1}{2}}(x;H)+\zeta(\tfrac{2}{3})\mathscr D_{\frac{1}{3}}(x;H)\right)\right|^2dx
\end{align*}
might be able to reduce the upper bound in Theorem \ref{thm:QRH-equiv-variance-sqfull}, by replacing $X$ by $H$. However, we have the series of results \cite{Ivic-divisor-short}, \cite{Lester} for the original divisor problem that states that the variance of the divisor problem
\begin{align*}
\dfrac{1}{X}\stretchint{7ex}_{\bs X}^{2X} \left|\Delta_{1,1}(x+H)-\Delta_{1,1}(x)\right|^2 dx:=\dfrac{1}{X}\stretchint{7ex}_{\bs X}^{2X} \left|\sum_{x<ab\leqslant x+H}1-\left(x\log x+(2\gamma-1)x\right)\right|^2 dx \gg H
\end{align*}
for $H\leqslant X^{\frac{1}{2}-\eps}$. In addition, Theorem 2.3 in the work of Lester \cite{Lester} showed that there exists a constant $B$ such that the left-hand side above is $\sim BX^{1/2}$ for $x^{1/2+\eps}\leqslant H\leqslant x^{1-\eps}$. This is peculiar, especially when $H=x^{1-\eps}$ since $\Delta_{1,1}(x+H)$ and $\Delta_{1,1}(x)$ seem uncorrelated. For this reason, we kept the bound to be a power of $X$ rather than of $H$.
\end{remark}

\begin{lemma}[Theorem 6.2 of \cite{Tenenbaum}]
\label{lem:intboundon_exp(2pi i)}
Let $f$ be a differentiable function in $(a,b)$ such that $f'(u)$ is monotonic of the constant sign and $m:=\inf_{u\in(a,b)}|f'(u)|>0$. Then,
$$\left|\stretchint{7ex}_{\bs a}^b e(f(u))du\right|\ll \dfrac{1}{m}.$$
\end{lemma}

\begin{lemma}
\label{lem:reducetoexp1}
Given positive integers $k,m$ such that $m>k$ and $\gcd(m,k)=1$. Let $c>0$ be a constant and $y\geqslant 1$. We obtain that
$$\sum_{d\leqslant y}{\vphantom{\sum}}' \left|c^{\frac{m}{k}}-d^{\frac{m}{k}}\right|^{-1}\ll_{k} \sum_{d\leqslant y}{\vphantom{\sum}}'|c-d|^{-1} \text{ and } \sum_{d\leqslant y}{\vphantom{\sum}}' |c^{\frac{2}{3}}-d^{\frac{2}{3}}|^{-1}\ll \sum_{d\leqslant y}{\vphantom{\sum}}'y^{\frac{1}{3}}|c-d|^{-1},$$
where the prime denotes the sum over $|c-d|> 1$.
\end{lemma}

\begin{proof}
For the first assertion, we see by completing the exponent to an integer that for the summand satisfying $d>c$
\begin{align*}
    \dfrac{1}{d^{\frac{m}{k}}-c^{\frac{m}{k}}}\leqslant \dfrac{kd^{m-\frac{m}{k}}}{(d-c)d^{m-1}}\ll_k \dfrac{1}{d-c}.
\end{align*}
The summand satisfying $d\leqslant c$ is done similarly by changing the role of $c$ and $d$. Thus, we obtain the desired bound. The second assertion can be done similarly, but since $\frac{2}{3}<1$, we have the extra term $y^{\frac{1}{3}}$ on the right-hand side.
\end{proof}

\begin{proof}[Proof of Proposition \ref{prop:var-(2,3)-divisor-problem}]

We have from the hyperbola method that the left-hand side of the integrand of Equation \eqref{eq:(2,3)-divisor problem-variance} can be written as
\begin{equation}
\label{eq:Delta23_fromHyperbola}
\displaystyle \sum_{d\leqslant x^{\frac{1}{5}}} \psi({\sqrt{x}/d^{\frac{3}{2}}})+\sum_{d\leqslant x^{\frac{1}{5}}} \psi({\sqrt[3]{x}/d^{\frac{2}{3}}})+O(1),
\end{equation}
where we write $\psi(u):=u-[u]-\tfrac{1}{2}$ for all $u\in \R$. It, thus, suffices to work on the sum of the types
\begin{equation}
\label{eq:doublesum32_23}
\sum_{d,d_1\leqslant X^{\frac{1}{5}}} \stretchint{7ex}_{\bs d'^5}^X \psi({\sqrt{x}/d^{\frac{3}{2}}})\psi({\sqrt{x}/d_1^{\frac{3}{2}}})dx \text{ and }
\sum_{d,d_1\leqslant X^{\frac{1}{5}}} \stretchint{7ex}_{\bs d'^5}^{X} \psi(\sqrt[3]{x}/d^{\frac{2}{3}})\psi({\sqrt[3]{x}/d_1^{\frac{2}{3}}})dx
\end{equation}
where $d'=\max(d,d_1)$. Here, it is convenient to write $|\widetilde{\psi}|(t):=|\widetilde{\psi}(t)|$ which is also a Lipschitz function by applying the triangle inequality. Also, we have that
\begin{equation}
\label{eq:approxB1byB}
\left|\psi\left(s\right)\psi\left(s_1\right)-\widetilde{\psi}\left(s\right)\widetilde{\psi}\left(s_1\right)\right| \leqslant h\left(s\right)h\left(s_1\right)+ |\widetilde{\psi}|\left(s\right)h\left(s_1\right)+h\left(s\right)|\widetilde{\psi}|\left(s_1\right),
\end{equation}
and we may replace $s=\sqrt{x}/d^{\frac{3}{2}}$ and $s_1=\sqrt{x}/d_1^{\frac{3}{2}}$ or $s=\sqrt[3]{x}/d^{\frac{2}{3}}$ and $s_1=\sqrt[3]{x}/d_1^{\frac{2}{3}}$ for the double sums in \eqref{eq:doublesum32_23}. Then, we find the Fourier series for the functions $\psi(x), h(x),|\widetilde{\psi}|(x)$ as follows.
\begin{align*}
\widetilde{\psi}(t) &= \sum_{j\geqslant 1} a_j\sin(2\pi jt), h(t) = \dfrac{1}{2J}+\sum_{j\geqslant 1}b_j\cos(2\pi jt), |\widetilde{\psi}|(t) = 1+\dfrac{1}{2J}+\sum_{j\geqslant 1} c_j\cos(2\pi jt)
\end{align*}
where $|a_j|,|b_j|,|c_j|\ll \min(j,J)/j^2$ for each $j\geqslant 1$. By Equation \eqref{eq:approxB1byB}, we may consider $\widetilde{\psi}\widetilde{\psi}, hh, |\widetilde{\psi}|h, h|\widetilde{\psi}|$ instead of $\psi\psi$, and each of them (for the first double sum in \eqref{eq:doublesum32_23}) can be bounded by
\begin{align}
\label{eq:main-term-3/2}
\ll \dfrac{X^{\frac{7}{5}}}{J}+\sum_{d,d_1\leqslant X^{\frac{1}{5}}} \sum_{j,j_1\geqslant 1} |a_jb_{j_1}|\cdot\sum_{\substack{k= O(\log x)\\ N=X/2^k}}\left|\stretchint{7ex}_{\bs N}^{2N} e(\sqrt{x}(j/d^{\frac{3}{2}}-j_1/d_1^{\frac{3}{2}}))dx\right|,
\end{align}
where $a_j, b_{j_1}$ are any terms from $a_j,b_j,c_j$ in the Fourier series above. Employing Lemma \ref{lem:intboundon_exp(2pi i)}, we have that
\begin{align}
\label{eq:exponent-bound-3/2}
\left|\stretchint{7ex}_{\bs N}^{2N} e(\sqrt{x}(j/d^{\frac{3}{2}}-j_1/d_1^{\frac{3}{2}}))dx\right|\ll \sqrt{N}d^{\frac{3}{2}}d_1^{\frac{3}{2}}|jd_1^{\frac{3}{2}}-j_1d^{\frac{3}{2}}|^{-1}.
\end{align}

We now focus on the following summation.
\begin{align*}
    \sideset{}{'}\sum_{d\leqslant X^{\frac{1}{5}}} |jd_1^{\frac{3}{2}}-j_1d^{\frac{3}{2}}|^{-1}=j_1^{-1} \sideset{}{'}\sum_{d\leqslant X^{\frac{1}{5}}} \left|d^{\frac{3}{2}}-(j/j_1)d_1^{\frac{3}{2}}\right|^{-1},
\end{align*}
where the sum is over $j\geqslant j_1$ and $|d-(\tfrac{j}{j_1})^{\frac{2}{3}}d_1|>1$. The sum over $|d-(\tfrac{j}{j_1})^{\frac{2}{3}}d_1|\leqslant 1$ can be done by counting $d$ just once since there is at most one $d$ satisfying the inequality as $j/j_1\geqslant 1$. In this case, the contribution to \eqref{eq:main-term-3/2} is $O(X^{\frac{6}{5}})$ since $d$ is dependent on $d_1$. The case when $j<j_1$ can be done similarly by changing the role of $d,d_1$.\\

We have the identity
$$ \left|d^{\frac{3}{2}}-(j/j_1)d_1^{\frac{3}{2}}\right|^{-1} =d^{-\frac{1}{2}+\frac{1}{m}}\left|d^{1+\frac{1}{m}}-\left((j/j_1)^{\frac{m}{m+1}}\frac{d_1^{\frac{3m}{2(m+1)}}}{d^{\frac{m-2}{2(m+1)}}}\right)^{1+\frac{1}{m}}\right|^{-1},$$
and by absorbing the term $d^{-1/2+1/m}$ into the right-hand side of \eqref{eq:exponent-bound-3/2} we can consider summing over the reciprocal of absolute value. Thus by Lemma \ref{lem:reducetoexp1}, we have for any $m\in\mathbb N$ large,
\begin{align*}
   \sideset{}{'}\sum_{d\leqslant X^{\frac{1}{5}}} \left|d^{1+\frac{1}{m}}-\left((j/j_1)^{\frac{m}{m+1}}\frac{d_1^{\frac{3m}{2(m+1)}}}{d^{\frac{m-2}{2(m+1)}}}\right)^{1+\frac{1}{m}}\right|^{-1}\ll_m \sideset{}{'}\sum_{d\leqslant X^{\frac{1}{5}}} \left|d-c\right|^{-1},
\end{align*}
where $c=(j/j_1)^{\frac{m}{m+1}}\frac{d_1^{\frac{3m}{2(m+1)}}}{d^{\frac{m-2}{2(m+1)}}}$. Recall that we consider when $|d-(j/j_1)^{\frac{2}{3}}d_1|>1$. Suppose that $d>(j/j_1)^{\frac{2}{3}}d_1$, we have that
\begin{align*}
    c=(j/j_1)^{\frac{m}{m+1}}\frac{d_1^{\frac{3m}{2(m+1)}}}{d^{\frac{m-2}{2(m+1)}}}<\left(\frac{j}{j_1}\right)^{\frac{2}{3}}d_1<d.
\end{align*}
Therefore, in this case we obtain
\begin{align*}
|c-d|^{-1}=(d-c)^{-1}\leqslant (d-(j/j_1)^{\frac{2}{3}}d_1)^{-1}=|d-(j/j_1)^{\frac{2}{3}}d_1|^{-1}.    
\end{align*}
When $d\leqslant (j/j_1)^{\frac{2}{3}}d_1$, we will obtain similarly with the opposite direction so that we have
\begin{align*}
   \sideset{}{'}\sum_{d\leqslant X^{\frac{1}{5}}} \left|d^{1+\frac{1}{m}}-\left((j/j_1)^{\frac{m}{m+1}}\frac{d_1^{\frac{3m}{2(m+1)}}}{d^{\frac{m-2}{2(m+1)}}}\right)^{1+\frac{1}{m}}\right|^{-1}\ll_m \sideset{}{'}\sum_{d\leqslant X^{\frac{1}{5}}} \left|d-(j/j_1)^{\frac{2}{3}}d_1\right|^{-1}=O_\eta(\log X),
\end{align*}
where $\eta=1/m$ is small. The first inequality above holds from replacing $d$ by $d_1$ by inspecting the case for the absolute value. Similar arguments hold for the other case in \eqref{eq:doublesum32_23}. The bound we obtained is uniform on $j,j_1$. Thus, the left-hand side of Equation \eqref{eq:exponent-bound-3/2} is at most
\begin{align*}
\ll_\eps (\log X)\sqrt{N}\cdot X^{\frac{1}{5}+\eps}\cdot \left(\sum_{d_1\leqslant X^{\frac{1}{5}}} d_1^{\frac{3}{2}} \right)\ll_\eps X^{\frac{6}{5}+\eps},
\end{align*}
since the sum is over $d\leqslant X^{\frac{1}{5}}$ and $N\leqslant X, d_1\leqslant X^{\frac{1}{5}}$. This proves the proposition.
\end{proof}

Finally, we prove that $\frac{1}{10}$-QRH with Proposition \ref{prop:var-(2,3)-divisor-problem} allows us to compute the expected order of the variance of square-full integers in short intervals. Let us first present lemmas used in computing the variance. The first two can be proved by adapting the proofs of Littlewood for analogous results involving the full RH; see, for example, \cite{Titchmarsh}. 
\begin{lemma}[Littlewood-type equivalence form of the $\frac{1}{10}$-QRH]
\label{lem:1/10-QRH-equiv-mobius-sum}
The $\frac{1}{10}$-QRH is equivalent to the following statement. For any $\eps>0$, we have that
$$\sum_{n\leqslant x}\mu(n)=O_\eps(x^{\frac{3}{5}+\eps}).$$
\end{lemma}

\begin{lemma}[Subconvexity bound under $\frac{1}{10}$-LH]
\label{lem:subconvexity-sigma>3/5}
Assuming the $\frac{1}{10}$-QRH, we have for $\sigma\geqslant\tfrac{3}{5}$ that
$$|\zeta(\sigma+it)|\ll_\eps (1+|t|)^\eps \text{ and }\dfrac{1}{|\zeta(\sigma+it)|}\ll_\eps (1+|t|)^\eps.$$
\end{lemma}

\begin{lemma}[Saffari-Vaughan's inequality, {\cite[Page~25]{Saffari-Vaughan}}]
\label{lem:Saffari-Vaughan}
Let $F:\R\to\mathbb C$ be a square-integrable function, and $H\leqslant X$. Then, we have that
$$\stretchint{7ex}_{\bs X}^{2X}|F(x+H)-F(x)|^2\ll \sup_{\theta\in[\frac{H}{3X}, \frac{3H}{X}]}\stretchint{7ex}_{\bs X}^{3X}|F(u(1+\theta))-F(u)|^2du.$$
\end{lemma}

\begin{remark}
We shall use Lemma \ref{lem:subconvexity-sigma>3/5} to prove the following statement. Let $\eps>0$. Suppose that $T\geqslant 1$ and $|t|\leqslant T$. Then, under the $\frac{1}{10}$-QRH we have that
\begin{equation}
\label{eq:summatory_mu(n)/n^it}
\sum_{n\leqslant x}\dfrac{\mu(n)}{n^{it}}\ll_\eps x^{\frac{3}{5}}(xT)^{\eps}.
\end{equation}
This can be seen as follows. Using Perron's formula, we obtain
\begin{align*}
\sum_{n\leqslant x} \dfrac{\mu(n)}{n^{it}}=\dfrac{1}{2\pi i}\stretchint{7ex}_{\bs 2-iW}^{2+iW} f(w)\dfrac{x^w}{w}dw+O(x^2W^{-1}+1),
\end{align*}
where $$f(w):=\sum_{n\geqslant 1}\dfrac{\mu(n)n^{-it}}{n^w}=\zeta(w+it)^{-1}.$$
Therefore, by choosing $W=x^2$ and moving the contour to the vertical line $\Re w=\tfrac{3}{5}$, we obtain the desired result by invoking Lemma \ref{lem:subconvexity-sigma>3/5}.
\end{remark}

\subsection{$\frac{1}{10}$-QRH implies the result on the variance}

First of all, we can write the left-hand side of Equation \eqref{eq:Variance-of-full-short-ranges} as
\begin{equation}
\label{term:var-short-interval-sum-on-a2b3c6}
\dfrac{1}{X}\stretchint{7ex}_{\bs X}^{2X}\left|\sum_{x<nc^6\leqslant x+H}\mu(c)1_{2,3}(n)-\left(\dfrac{\zeta(\frac{3}{2})}{\zeta(3)}(\sqrt{x+H}-\sqrt{x})+\dfrac{\zeta(\frac{2}{3})}{\zeta(2)}(\sqrt[3]{x+H}-\sqrt[3]{x})\right)\right|^2dx.
\end{equation}
It is natural to define
\begin{align}
\label{eq:A(x)-for-C<H^1/6}
A(x):= \sum_{\substack{nc^6\leqslant x\\ C<c\leqslant 2C}}\mu(c)1_{2,3}(n)-\left(\sum_{C<c\leqslant 2C}\dfrac{\zeta(\frac{3}{2})\mu(c)}{c^3}x^{\frac{1}{2}}+\sum_{C<c\leqslant 2C}\dfrac{\zeta(\frac{2}{3})\mu(c)}{c^2}x^{\frac{1}{3}}\right),
\end{align}
where the possibility of $C$ is split into cases: $C\leqslant H^{\frac{1}{6}}$ and $C\geqslant H^{\frac{1}{6}}$.\\

\noindent
\textbf{Case $C\leqslant H^{\frac{1}{6}}$.}
\noindent
We split the case further according to the size of $|t|$ into $3$ cases, namely $|t|\leqslant (X/C^{6})H^{-\frac{\eps}{100}}$ and $|t|\geqslant (X/C^{6})H^{-\frac{\eps}{100}}$.\\

\noindent
1) $|t|\leqslant (X/C^{6})H^{-\frac{\eps}{100}}$. In this subcase, we shall modify $A(x)$ in Equation \eqref{eq:A(x)-for-C<H^1/6} slightly. Let $f$ be a smooth function supported on $[\frac{1}{20},20]$ such that $f\equiv 1$ on $[\frac{1}{10},10]$ and its derivative $f'\in L^{1}(\R)$. Notice that we can plug in this function into $A(x)$ in \eqref{eq:A(x)-for-C<H^1/6} since the initial terms cancel when we sum over the range $[x,x+H]$. Define
\begin{align}
\label{eq:A(x)-forC<H1/6-proof}
A(x):= \sum_{\substack{nc^6\leqslant x\\ C<c\leqslant 2C\\ n=a^2b^3}}\mu(c)f\left(\dfrac{n}{X/C^6}\right)-\left(\sum_{C<c\leqslant 2C}\dfrac{\zeta(\frac{3}{2})\mu(c)}{c^3}x^{\frac{1}{2}}+\sum_{C<c\leqslant 2C}\dfrac{\zeta(\frac{2}{3})\mu(c)}{c^2}x^{\frac{1}{3}}\right).
\end{align}

By Perron's formula, we obtain
\begin{equation}
\label{eq:A(x)-C<H1/6-proof-perron}
A(x)=\stretchint{7ex}_{\bs 2-i\infty}^{2+i\infty} \dfrac{x^s}{s}\cdot N_1(s)M(6s)ds-\left(\sum_{C<c\leqslant 2C}\dfrac{\zeta(\frac{3}{2})\mu(c)}{c^3}x^{\frac{1}{2}}+\sum_{C<c\leqslant 2C}\dfrac{\zeta(\frac{2}{3})\mu(c)}{c^2}x^{\frac{1}{3}}\right),
\end{equation}
where 
\begin{align}
\label{eq:M-N1-defined}
M(s):=\sum_{C<c\leqslant 2C} \dfrac{\mu(c)}{c^s} \text{ and } N_1(s) :=\sum_{\substack{n=a^2b^3}}\dfrac{1}{n^s}f\left(\dfrac{n}{X/C^6}\right).
\end{align}

Now, we state some useful functions for this special subcase of $t$. Write
\begin{align}
\label{eq:N1tilde-N2}
\widetilde{N}_1(s) &:=\sum_{\substack{n}} \left(\dfrac{\zeta(\frac{3}{2})}{2n^{s+\frac{1}{2}}}+\dfrac{\zeta(\frac{2}{3})}{3n^{s+\frac{2}{3}}}\right)f\left(\dfrac{n}{X/C^6}\right) \nonumber
\\ N_2(s) &:=\stretchint{7ex}_{\bs\mathbb R}\left(\dfrac{\zeta(\frac{3}{2})}{2u^{s+\frac{1}{2}}}+\dfrac{\zeta(\frac{2}{3})}{3u^{s+\frac{2}{3}}}\right)f\left(\dfrac{u}{X/C^6}\right)du,
\end{align}
so that we can use $\widetilde{N}_1(s)$ as an approximation of $N_1(s)$. This function is smooth on the domain of the integral, for which we can employ Poisson summation. Also, note that we have that $N_1-N_2=(N_1-\widetilde{N}_1)+(\widetilde{N}_1-N_2)=G+(\widetilde{N}_1-N_2)$, where 
\begin{equation}
\label{eq:G-inproof-C<H1/6}
G(s) := \sum_{n}\dfrac{g(n)}{n^s}f\left(\dfrac{n}{X/C^6}\right),
\end{equation}
and $g(n):=1_{2,3}(n)-\frac{\zeta(\frac{3}{2})}{2n^{\frac{1}{2}}}-\frac{\zeta(\frac{2}{3})}{3n^{\frac{2}{3}}}$. We see that
\begin{align*}
\stretchint{7ex}_{\bs\frac{1}{2}-i\infty}^{\frac{1}{2}+i\infty} \dfrac{x^s}{s}N_2(s)M(6s)ds=x^\frac{1}{2}\sum_{C<c\leqslant 2C} \dfrac{\zeta(\frac{3}{2})\mu(c)}{c^{3}} +x^{\frac{1}{3}}\sum_{C<c\leqslant 2C} \dfrac{\zeta(\frac{2}{3})\mu(c)}{c^{2}}.
\end{align*}
Therefore, upon changing the path in the integral of $A(x)$ in Equation \eqref{eq:A(x)-C<H1/6-proof-perron} to $\frac{1}{2}+it$ and plugging in the above equation, we have that
\begin{equation}
\label{eq:A(x)-contour-caseC<H1/6-t<X/C6}
A(x)=\stretchint{7ex}_{\bs\R}\dfrac{x^{\frac{1}{2}+it}}{\frac{1}{2}+it}\left(N_1(\tfrac{1}{2}+it)-N_2(\tfrac{1}{2}+it)\right)\cdot M(3+6it)dt.
\end{equation}
We obtain by the Saffari-Vaughan's inequality (Lemma \ref{lem:Saffari-Vaughan}) that
\begin{align*}
\dfrac{1}{X}\stretchint{7ex}_{\bs X}^{2X} |A(x+H)-A(x)|^2 dx \ll \dfrac{1}{X}\stretchint{7ex}_{\bs X}^{3X} |A(u(1+\theta))-A(u)|^2 du,
\end{align*}
for some $\theta\in[\tfrac{H}{3X},\tfrac{3H}{X}]$. Now, let $u=e^{x}$ and $w$ be such that $1+\theta=e^w$, and thus, $w\asymp \tfrac{H}{X}$. By Equation \eqref{eq:A(x)-contour-caseC<H1/6-t<X/C6} with the above and changing variable, we have that
\begin{align}
\label{eq:bound-on-(N1-N2)M}
\dfrac{1}{X}&\stretchint{7ex}_{\bs X}^{2X}|A(x+H)-A(x)|^2 dx \nonumber\\
&\leqslant \dfrac{1}{X}\stretchint{7ex}_{\bs \log X}^{\log (3X)} e^{2x}\stretchint{7ex}_{\bs \R} \left|\dfrac{e^{w(\frac{1}{2}+it)}-1}{\frac{1}{2}+it}\right|^2 \left|\big(N_1(\tfrac{1}{2}+it)-N_2(\tfrac{1}{2}+it)\big)M(3+6it)\right|^2dt dx \nonumber\\
&\ll X\stretchint{7ex}_{\bs \mathbb{R}}\left|\dfrac{e^{w(\frac{1}{2}+it)}-1}{\frac{1}{2}+it}\right|^2\cdot \left|\big(N_1(\tfrac{1}{2}+it)-N_2(\tfrac{1}{2}+it)\big)M(3+6it)\right|^2dt.
\end{align}

We claim that $|\widetilde{N}_1(\tfrac{1}{2}+it)-N_2(\tfrac{1}{2}+it)|\ll_A X^{-A}$ for any $A>0$ using the Poisson summation formula. Hence, this part is negligible. By invoking the Cauchy-Schwarz inequality, we are left to prove the term involving $|G(\tfrac{1}{2}+it)|^2$. Now, we prove the former claim as follows. Write $\widetilde N_1(s):=N_{1,1}(s)+N_{1,2}(s)$ and $N_2(s):=N_{2,1}(s)+N_{2,2}(s)$ for the terms involving $n^{\frac{1}{2}}$ and $n^\frac{1}{3}$, respectively. Thus,
\begin{align*}
&N_{1,1} (\tfrac{1}{2}+it)=\sum_{m} \dfrac{\zeta(\frac{3}{2})}{2m^{1+it}}f\left(\dfrac{m}{X/C^6}\right)
\\ &=\sum_{\ell}\stretchint{7ex}_{\bs \R}\dfrac{\zeta(\frac{3}{2})}{2u^{1+it}}f\left(\dfrac{u}{X/C^6}\right)e(\ell u)du =N_{2,1}(\tfrac{1}{2}+it)+\sum_{\ell\not =0}\stretchint{7ex}_{\bs\R} \dfrac{\zeta(\frac{3}{2})}{2u^{1+it}}f\left(\dfrac{u}{X/C^6}\right)e(\ell u)du.
\end{align*}
To show $|N_{1,1}(1/2+it)-N_{2,1}(1/2+it)|\ll_{A,\eps} X^{-A}$, it suffices to show the same upper bound on
\begin{align*}
\dfrac{X}{C^6}\sum_{\ell\not=0}\stretchint{7ex}_{\bs 1/100}^{100} \left(\dfrac{C^6}{X}\right)^{1+it}\cdot\dfrac{f(y)}{y^{1+it}}e\left(\dfrac{\ell yX}{C^6}\right) dy,
\end{align*}
which is possible using partial integration on $e(\ell yX/C^6)$ repeatedly, and note that $|t|\leqslant (X/C^{6})H^{-\frac{\eps}{100}}$. This constraint is used since the term $t$ will appear in the numerator several times depending on $A$. Remember that $H\geqslant X^\eps$, we can achieve the upper bound. Similarly, we can obtain the same bound for $|N_{1,2}(\tfrac{1}{2}+it)-N_{2,2}(\tfrac{1}{2}+it)|$.\\

For the term involving $|G(\tfrac{1}{2}+it)|^2$, we recall that we want to determine the term left from Equation \eqref{eq:bound-on-(N1-N2)M}
\begin{align*}
X\stretchint{7ex}_{\bs \mathbb{R}}\left|\dfrac{e^{w(\frac{1}{2}+it)}-1}{\frac{1}{2}+it}\right|^2\cdot |G(\tfrac{1}{2}+it)|^2|M(3+6it)|^2 dx,
\end{align*}
where we can apply partial summation on Equation \eqref{eq:summatory_mu(n)/n^it} involving the sum on the M\"{o}bius function to the last term above to obtain the bound $\ll_\eps C^{-\frac{24}{5}}X^{\eps}$ since $|t|\leqslant X$. It suffices to establish the remaining term to be $\ll_\eps X^{\frac{1}{5}+\eps}C^{\frac{24}{5}}$. Now, define
\begin{equation}
B(x):=\sum_{n\leqslant x} g(n)f\left(\frac{n}{X/C^6}\right)=\sum_{n\leqslant x} \left(1_{2,3}(n)-\left(\frac{\zeta(\tfrac{3}{2})}{2n^{\frac{1}{2}}}+\frac{\zeta(\tfrac{2}{3})}{3n^{\frac{2}{3}}}\right)\right)f\left(\dfrac{n}{X/C^6}\right),
\end{equation}
we have by Perron's formula that
\begin{align}
\label{eq:Perron-B}
&\dfrac{B(e^{w+x})-B(e^x)}{e^{x/2}}=\dfrac{1}{2\pi}\stretchint{7ex}_{\bs\mathbb{R}} \left(\dfrac{e^{w(\frac{1}{2}+it)}-1}{\frac{1}{2}+it}\right)e^{itx}G(\tfrac{1}{2}+it) dt.
\end{align}
Let us now write 
\begin{align*}
    h(t):= \left(\dfrac{e^{w(\frac{1}{2}+it)}-1}{\frac{1}{2}+it}\right)G(\tfrac{1}{2}+it).
\end{align*}
Then, we have that
\begin{align*}
    \stretchint{7ex}_{\bs \R} \left(\dfrac{e^{w(\frac{1}{2}+it)}-1}{\frac{1}{2}+it}\right)e^{itx}G(\tfrac{1}{2}+it)dt=\widehat{h}(-\frac{x}{2\pi})
\end{align*}
and by Plancherel's theorem
\begin{align*}
    \stretchint{7ex}_{\bs \R} |h(t)|^2 dt=\stretchint{7ex}_{\bs \R} |\widehat{h}(x)|^2 dx=\frac{1}{2\pi} \stretchint{7ex}_{\bs \R} \left| \widehat{h}(-\frac{x}{2\pi})\right|^2 dx=2\pi \stretchint{7ex}_{\bs \R} \frac{|B(e^{w+x})-B(e^x)|^2}{e^x} dx,
\end{align*}
where the last equality comes from integrating over $x$ for both sides of \eqref{eq:Perron-B}. Similar to before, we write $u=e^x$ and $1+\theta=e^w$ so that $w\asymp \tfrac{H}{X}$. Using the above, we obtain
\begin{align}
\label{eq:C<H1/6-x-(2,3)}
&\dfrac{X}{2\pi}\stretchint{7ex}_{\bs\mathbb{R}} \left|\dfrac{e^{w(\frac{1}{2}+it)}-1}{\tfrac{1}{2}+it}\right|^2|G(\tfrac{1}{2}+it)|^2 dt=X\stretchint{7ex}_{\bs 0}^\infty \dfrac{|B(e^{w+x})-B(e^x)|^2}{e^x}dx \nonumber\\
&\quad\ll X\stretchint{7ex}_{\bs 0}^\infty |B(u(1+\theta))-B(u)|^2 \dfrac{du}{u^2} \nonumber\\
&\quad\ll X\stretchint{7ex}_{\bs\frac{X}{10C^6}}^{\frac{10X}{C^6}}\dfrac{B^2(x)}{X^2/C^{12}} dx+X\stretchint{7ex}_{\bs \frac{X}{20C^6}}^{\frac{X}{10C^6}} \dfrac{1}{X^2/C^{12}} |B(x+H)-B(x)|^2 dx+J,
\end{align}
where $H=H(x)=o(x)$, and $J$ is the integral with the same integrand as in the second term but over $[\tfrac{10X}{C^6},\tfrac{20X}{C^6}]$. We now obtain
\begin{align*}
\dfrac{1}{X}\stretchint{7ex}_{\bs \tfrac{X}{20C^6}}^{\tfrac{X}{10C^6}}\left|B(x+H)-B(x)\right|^2 dx &= \dfrac{1}{X}\stretchint{7ex}_{\bs \tfrac{X}{20C^6}}^{\tfrac{X}{10C^6}} \left|\sum_{x<n\leqslant x+H}1_{2,3}(n)-\left(\frac{\zeta(\tfrac{3}{2})}{2n^{\frac{1}{2}}}+\frac{\zeta(\tfrac{2}{3})}{3n^{\frac{2}{3}}}\right)\right|^2 dx\\
&\ll \dfrac{1}{X}\stretchint{7ex}_{\bs\tfrac{X}{20C^6}}^{\tfrac{X}{10C^6}} \big(\Delta^2(x)+\Delta^2(x+H)+O(1)\big)dx,
\end{align*}
and the term $J$ can be done similarly. Thus, Equation \eqref{eq:C<H1/6-x-(2,3)} is bounded above by $O_\eps(X^{\frac{1}{5}+\eps}C^{\frac{24}{5}})$ by invoking Proposition \ref{prop:var-(2,3)-divisor-problem} as desired.\\

\noindent
2) $|t|\geqslant (X/C^{6})H^{-\frac{\eps}{100}}$. We have by the definition in \eqref{eq:M-N1-defined} that 
\begin{align*}
    N_1(\tfrac{1}{2}+it)\ll \sum_{\frac{X}{20C^6}<a^2b^3\leqslant \frac{20X}{C^6}} \dfrac{1}{\sqrt{n}}\ll \dfrac{1}{\sqrt{X/C^6}}\cdot Q(20X/C^6)\ll 1.
\end{align*}
Here, we recall that $Q(x)$ is the number of square-full integers not exceeding $x$. Also, $N_2(\frac{1}{2}+it)$ can be done similarly, and the similar bound is obtained. Now, we restrict to the subcase when $|t|\leqslant X^{10}$. By $\frac{1}{10}$-QRH, we also have from performing partial summation as in \eqref{eq:summatory_mu(n)/n^it} to obtain $|M(3+6it)|\ll_\eps C^{-\frac{12}{5}+\frac{\eps}{2}}$. Note that, in this case, $1/|t|\leqslant C^{6}H^\frac{\eps}{100}/X$. Hence, the contribution of this case to the integral on the right-hand side of \eqref{eq:bound-on-(N1-N2)M} is 
\begin{align*}
    \ll_\eps &XC^{-\frac{24}{5}+\eps}\stretchint{7ex}_{\bs (X/C^6)H^{-\frac{\eps}{100}}\leqslant |t|\leqslant X^{10}}\frac{1}{|t|^2} dt  \ll C^{\frac{6}{5}+\eps}H^{\frac{\eps}{100}}\stretchint{7ex}_{|t|\leqslant X^{10}} \frac{dt}{|t|}\ll_\eps H^{\frac{1}{5}+\eps},
\end{align*}
which is the desired bound as $C\leqslant H^{\frac{1}{6}}$ and $\log X\ll_\eps X^{\eps^2/6} \ll_\eps H^{\eps/6}$. When $|t|\geqslant X^{10}$, we use the similar bound for $N_1(\tfrac{1}{2}+it),N_2(\tfrac{1}{2}+it)$ and $|M(3+6it)|\ll 1$. Thus, the right-hand side of \eqref{eq:bound-on-(N1-N2)M} for this case is $O(1)$.\\

\noindent
\textbf{Case $C\geqslant H^{\frac{1}{6}}$.} Similar to defining $A(x)$ as in \eqref{eq:A(x)-for-C<H^1/6}, we define the following tail for the integrand in \eqref{term:var-short-interval-sum-on-a2b3c6} by
\begin{align*}
    A(x):=\sum_{\substack{a^2b^3c^6\leqslant x\\ c>C}}\mu(c)-\left(x^{\frac{1}{2}}\sum_{c>C} \dfrac{\zeta(\frac{3}{2})\mu(c)}{c^{3}}+x^{\frac{1}{3}}\sum_{c>C} \dfrac{\zeta(\frac{2}{3})\mu(c)}{c^{2}}\right),
\end{align*}
as we now do not need to use Poisson summation to bound the difference between the Dirichlet polynomials $N_1,N_2$ by $X^{-A}$, we replace $f$ in the previous case simply by $1$. We obtain via Perron's formula 
\begin{align*}
    \dfrac{A(e^{w+x})-A(e^x)}{e^{x/2}}=\dfrac{1}{2\pi}\stretchint{7ex}_{\bs\R} \left(\dfrac{e^{w(\frac{1}{2}+it)}-1}{\frac{1}{2}+it}\right)e^{itx}\zeta(1+2it)\zeta(\tfrac{3}{2}+3it)M(3+6it) dt,
\end{align*}
where $M(s)$ is similar to $M(s)$ in \eqref{eq:M-N1-defined}, which is now defined by
$$M(s):=\sum_{c>C} \dfrac{\mu(c)}{c^s}.$$
Then, write $e^{w}=1+\theta$ with $w\asymp \tfrac{H}{X}$ and use Saffari-Vaughan's inequality (Lemma \ref{lem:Saffari-Vaughan}) and Plancherel's identity similar to \eqref{eq:bound-on-(N1-N2)M} and the first equation in \eqref{eq:C<H1/6-x-(2,3)} to obtain
\begin{align*}
\dfrac{1}{X}&\stretchint{7ex}_{\bs X}^{2X}|A(x+H)-A(x)|^2 dx \ll X\stretchint{7ex}_{\bs 0}^\infty |A(u(1+\theta))-A(u)|^2\dfrac{du}{u^2} \nonumber \\
&\qquad \ll X\stretchint{7ex}_{\bs\R}
\min\left\{\left(\dfrac{H}{X}\right)^2,\dfrac{1}{|t|^2}\right\}|\zeta(1+2it)\zeta(\tfrac{3}{2}+3it)M(3+6it)|^2 dt\nonumber \\
&\qquad \ll \frac{H^2}{X}\stretchint{7ex}_{\bs |t|\leqslant X/H}|\zeta(1+2it)M(3+6it)|^2dt+X\stretchint{7ex}_{\bs \frac{X}{H}\leqslant |t|\leqslant X^2}|\zeta(1+2it)M(3+6it)|^2dt+1.
\end{align*}

Assuming $\frac{1}{10}$-QRH, we have from \eqref{eq:summatory_mu(n)/n^it} that $|M(3+6it)|\ll_\eps C^{-\frac{12}{5}}H^{\frac{\eps}{10}}$ for $|t|\leqslant X^{10}$, and the above implies 
\begin{align*}
\dfrac{1}{X}&\stretchint{7ex}_{\bs X}^{2X}|A(x+H)-A(x)|^2 dx \ll_\eps \dfrac{H^{1+\frac{\eps}{10}}\log^{10}X}{C^{\frac{24}{5}}}+1\ll H^{\frac{1}{5}+\eps},
\end{align*}
for $C\geqslant H^{1/6}$. The first integral employs variance of the zeta function, and the second term comes from the following arguments. We see that
\begin{align*}
    X &\stretchint{7ex}_{\bs\frac{X}{H}\leqslant |t|\leqslant X^2} |\zeta(1+2it)M(3+6it)|^2dt\\
    &\quad \ll X\stretchint{7ex}_{\bs \frac{X}{H}}^{X^2} \frac{1}{T} \stretchint{7ex}_{\bs T\leqslant |t|\leqslant 2T} \dfrac{1}{T^2}\cdot |\zeta(1+2it)M(3+6it)|^2dt\, dT\\
    &\quad \leqslant H\log (XH)\sup_{\frac{H}{X}\leqslant T\leqslant X^2} \dfrac{1}{T}\stretchint{7ex}_{\bs T\leqslant |t|\leqslant 2T} |\zeta(1+2it)M(3+6it)|^2dt\ll C^{-\frac{24}{5}}H^{\frac{\eps}{10}}\cdot H\log^{10}X,
\end{align*}
by the variance of the zeta function on the line $\Re s=1$. This completes the proof of the implication, thereby concluding the proof of Theorem \ref{thm:QRH-equiv-variance-sqfull}.

\section{\texorpdfstring{$q$}{q}-aspect and Variance of square-full integers in arithmetic progressions}
\label{sec:Propositions in the $q$-aspect and Variance of square-full integers in arithmetic progressions}

In this section, we provide and prove several propositions in the $q$-aspect, which are used in proving Theorem \ref{thm:sqfull-in-arith}. Since we are working with the nonlinear nature of $1^\blacksquare(n)$, we approximate its relative function $1_{2,3}(n)$ by a smooth function. This function $1_{2,3}(n)$ returns $1$ if $n$ can be written as $a^2b^3$ for some integers $a,b$ and zero otherwise. This is approximable in a long range and corresponds to the $(2,3)$-divisor problem as discussed in the previous section. We write for $n\in\N$
\begin{align}
\label{eq:approx=characteristic-sqfull}
1_{2,3}(n)=\left( \dfrac{\zeta(\tfrac{3}{2})}{2n^{\frac{1}{2}}}+\dfrac{\zeta(\tfrac{2}{3})}{3n^{\frac{2}{3}}}\right)+\left(1_{2,3}(n)-\left( \dfrac{\zeta(\tfrac{3}{2})}{2n^{\frac{1}{2}}}+\dfrac{\zeta(\tfrac{2}{3})}{3n^{\frac{2}{3}}}\right)\right)=:f(n)+r(n).
\end{align}

The error that occurred from the approximation is handled by employing the unconditional bound of the $(2,3)$-divisor problem, which is the result involving the residual function $r(n)$, i.e.,
\begin{equation}
\label{eq:err-23-div}
    \sum_{n\leqslant x} r(n)= O(x^{\frac{2}{15}}).
\end{equation}

Let $b_n$ be a sequence of complex numbers, and $\alpha\in(\Z/q\Z)^\times$. We have the following identities that will be used throughout this section, which are proved by the orthogonality of characters.

\begin{align}
\label{eq:identities-q-aspect}
    \dfrac{1}{\varphi(q)}\sum_{\substack{\chi\not=\chi_0}} \dfrac{1+\chi(\alpha)}{2}\left| \sum_{n} b_n\chi(n)\right|^2 &= \sum_{\ell\in \qZ}\left| \frac{1}{2}\sum_{n\equiv\ell,\alpha\ell(q)} b_n-\dfrac{1}{\varphi(q)}\sum_{\gcd(n,q)=1}b_n\right|^2 \text{ and}\\
    \label{eq:identoties-q-aspect2}
    \dfrac{1}{\varphi(q)}\sum_{\substack{\chi(q)}}\frac{1+\chi(\alpha)}{2} \left| \sum_{n} b_n\chi(n)\right|^2 &= \sum_{\ell\in \qZ}\left|\frac{1}{2}\sum_{n\equiv\ell,\alpha\ell(q)} b_n\right|^2.
\end{align}

In proving Theorem \ref{thm:sqfull-in-arith}, we split the summand into pieces using approximation in Equation \eqref{eq:approx=characteristic-sqfull} involving $1_{2,3}(n)=f(n)+r(n)$. Note that the identity $\mu^2(n)=\sum_{c^2|n} \mu(c)$ implies that

\begin{align}
    \label{eq:replace-sqfull-by1_23}
    \sum_{\substack{n\equiv\ell(q)\\ x<n\leqslant 2x}}1^\blacksquare(n)=\sum_{\substack{nc^6\equiv\ell(q)\\ x<nc^6\leqslant 2x}}\mu(c)1_{2,3}(n).
\end{align}

Then, we split the variance summand into summing over certain thresholds on $c$, see Proposition \ref{prop:equiv-to-gcd-smallc} and Proposition \ref{prop:equiv-to-gcd-largec} below. These results cover each range of $c$ over the smooth main term $f(n)$.\\

In the proof of Proposition \ref{prop:equiv-to-gcd-smallc}, we employ the P\'olya formula to convert the character sum over integers in the interval into the sum over blocks of size $q$. Also, we use Poisson summation to exploit the Fourier transform to obtain a bound on the parameter $c$ as in \eqref{eq:replace-sqfull-by1_23} from below. Then, we use identities at the beginning and the majorant principle to replace involved functions by their majorant without disturbing the characters. In certain cases, we count the number of solutions occurring from the condition $nc^6\equiv \ell \pmod{q}$ from non-diagonal terms. For Proposition \ref{prop:equiv-to-gcd-largec}, we use the Perron formula and bound the Dirichlet polynomial using the large-value theorem (Lemma \ref{lem:Hybrid large-value theorem} below) on a well-spaced set. 

\begin{proposition}
\label{prop:equiv-to-gcd-smallc}
Let $\eps\in(0,\frac{1}{1000})$. Let $x^{\frac{1}{11}}\leqslant q\leqslant x^{1-50\eps}$ be a prime and $(x/q)^{1+\eps}\leqslant z\leqslant x^{-\eps}\sqrt{qx}$. For any quadratic nonresidue $\alpha$ modulo $q$, we have that
\begin{align}
\label{eq:prop2.2}
    \dfrac{1}{\varphi(q)}\sum_{\ell\in(\Z/q\Z)^\times}\left|\dfrac{1}{2}\sum_{\substack{x<nc^6\leqslant 2x\\ nc^6\equiv \ell,\alpha\ell(q)\\ c^6\leqslant z}} \mu(c)f(n)-\frac{1}{\varphi(q)}\sum_{\substack{x<nc^6\leqslant 2x\\\gcd(nc^6,q)=1\\ c^6\leqslant z}}\mu(c)f(n)\right|^2\sim C(x/q)^{\frac{1}{6}},
\end{align}
where $C$ is defined as in Theorem \ref{thm:sqfull-in-arith} and $f$ is as in \eqref{eq:approx=characteristic-sqfull}.
\end{proposition}

\begin{proposition}
\label{prop:equiv-to-gcd-largec}
Let $\eps\in(0,\frac{1}{1000})$ and $q$ be a prime such that $q\gg x^{1/9+\eps}$ and $z\geqslant (x/q)^{1+\eps}$, then for $f$ and $\alpha$ as in Proposition \ref{prop:equiv-to-gcd-smallc}
\begin{align*}
    \dfrac{1}{\varphi(q)}\sum_{\ell\in(\Z/q\Z)^\times}\left|\frac{1}{2}\sum_{\substack{x<nc^6\leqslant 2x\\ nc^6\equiv \ell,\alpha\ell(q)\\ c^6> z}} \mu(c)f(n)-\dfrac{1}{\varphi(q)}\sum_{\substack{x<nc^6\leqslant 2x\\\gcd(nc^6,q)=1\\ c^6> z}}\mu(c)f(n)\right|^2\ll_\eps (x/q)^{\frac{1}{6}-\frac{\eps}{16}}.
\end{align*}
\end{proposition}

Let us prove each proposition. First, we state lemmas involving the subconvexity bound, the estimate of the fourth moment average of $L$-functions, and the size of the well-spaced set, i.e., the set of tuples $(t,\chi), (t',\chi')$ such that $|t-t'|\geqslant 1$ or $\chi\not=\chi'$, for which the Dirichlet polynomial is large, as follows.

\begin{lemma}[Hybrid subconvexity bound, {\cite[Theorem~1.1]{Petrow-Young}}]
\label{lem:hyb-subconvexity-sigma>1/2}
For a prime $q$, a primitive character $\chi\, (\text{mod }q)$, and $\sigma\geqslant\tfrac{1}{2}$ we have that
$$|L(\sigma+it,\chi)|\ll_\eps (1+(q|t|))^{\frac{1}{6}+\eps}.$$
\end{lemma}

\begin{lemma}[Hybrid fourth moment estimate, {\cite[Theorem~10.1]{Montgomery}}]
\label{lem:hyb-fourth-moment-subconvexity-sigma>1/2}
Let $T,q\geqslant 2$ and $\sigma\geqslant \frac{1}{2}$,
$$\sum_{\chi(q)}\stretchint{7ex}_{\bs |t|\leqslant T}| L(\sigma+it,\chi)|^4dt\ll T\varphi(q)\log^4(Tq).$$
\end{lemma}

\begin{lemma}[Hybrid large-value theorem]
\label{lem:Hybrid large-value theorem}
    Let $q\in\N, N,T\geqslant 1$ and $V>0$. Let $F(s,\chi)=\sum_{n\leqslant N} a_n\chi(n)n^{-s}$ be a Dirichlet polynomial, and $G=\sum_{n\leqslant N} |a_n|^2$. Let $\mathcal T$ be a well-spaced set of tuples $(t_r,\chi)$ with $|t_r|\leqslant T$ and $\chi$ a primitive character modulo $q$ such that $|F(it_r,\chi)|\geqslant V$. Then,
    $$|\mathcal T|\ll GNV^{-2}+qT\min(GV^{-2}, G^3NV^{-6})\log^{18}(2qNT).$$
\end{lemma}
\begin{proof}
    This follows from {\cite[Theorems 9.16 and 9.18 with $k=q$ and $Q=1$]{Iwaniec-Kowalski}}.
\end{proof}

\subsection{Proof of Proposition \ref{prop:equiv-to-gcd-smallc}}

In proving Proposition \ref{prop:equiv-to-gcd-smallc} we require Proposition \ref{prop:2.2.1} and Proposition \ref{prop:2.2.2} below, similarly as in the proof of Proposition $3$ of \cite{GoroMato}. The idea is to use Proposition \ref{prop:2.2.1}, which covers the sum over a specific small interval, to fill the gap from replacing the smooth function $g(x)$, which approximates the characteristic function of the whole relevant interval, back by such a characteristic function in Proposition \ref{prop:2.2.2}.\\

In fact, Proposition \ref{prop:2.2.1} is slightly stronger than what we will be using on average over a quadratic residue and nonresidue modulo $q$. By the identity \eqref{eq:identities-q-aspect} mentioned at the beginning of this section, we transform the left-hand side of Equation \eqref{eq:prop2.2} by the sum over characters, letting $b_m:=b^{(z)}_m:=\sum_{c^6\leqslant z, c^6|m} \mu(c)f(m/c^6)= \sum_{m=nc^6,c^6\leqslant z} \mu(c)f(n)$.

\begin{proposition}
    \label{prop:2.2.1}
    Let $q$ be a prime. For $q\leqslant x^{1-\eps}$, $z\leqslant 2x$ and $I\subset [x,3x]$,
    \begin{align*}
        \dfrac{1}{\varphi(q)}\sum_{\chi\not=\chi_0} \left| \sum_{\substack{c^6\leqslant z\\ nc^6\in I}}\mu(c)\chi(n)\chi(c)^6f(n)\right|^2\ll_\eps\left(z^{\frac{1}{3}}+|I|^{\frac{1}{6}}q^{\frac{5}{6}}\right)\log^6 x.
    \end{align*}
\end{proposition}

\begin{proof}
We use the P\'olya formula to split the interval into groups of size $q$ $$\sum_{n\in I/c^6}\chi(n)=\dfrac{\tau(\chi)}{2\pi i}\sum_{|n|\leqslant q} \overline\chi(n)g_{I/c^6}(n)+O(\log q),$$
where we write $I=[a,b]$ and 
\begin{align*}
    g_{I/c^6}(n) &:= \frac{1}{n}\left(e\left(\frac{na}{c^6q}\right)-e\left(\frac{nb}{c^6q}\right)\right).\\
    &=\frac{1}{n}e\left(\frac{nb}{c^6q}\right)\cdot \left(e\left(\frac{n(a-b)}{c^6q}\right)-1\right) \ll \left|\frac{1}{n}\right| \left|e\left(\frac{n(a-b)}{c^6q}\right)-1\right|.
\end{align*}
The above is $\ll 1/|n|$ and if $|n|\leqslant 10(c^6q/|I|)$, it is
$$
\ll \left| \frac{2\pi i(a-b)}{c^6q\cdot 1!}+\left(\frac{2\pi i (a-b)}{c^6q}\right)^2\frac{n}{2!}+\ldots\right|\ll |I|c^{-6}q^{-1},
$$
since from the second term onward, the terms are bounded by a constant multiple of the first term within this condition. We, therefore, define a majorant of $g_{I/c^6}(n)$ by $h_{I/c^6}(n)$ which takes the value $|I|c^{-6}q^{-1}$ if $|n|\leqslant 10(c^6q/|I|)$ and $1/n$ otherwise. Similarly, we write 
$$\displaystyle {g}_{I/c^6}(n;t):=\frac{1}{n}\left(e\left(\dfrac{na}{c^6q}\right)-e\left(\dfrac{nt}{q}\right)\right).$$ Then, we have by partial summation that
\begin{align*}
    \sum_{n\in I/c^6} \dfrac{\chi(n)}{\sqrt{n}}=c^3\left(\frac{1}{\sqrt{b}}-\frac{1}{\sqrt{a}}\right)\sum_{n\in I/c^6}\chi(n)+\dfrac{\tau(\chi)}{4\pi i}\sum_{|n|\leqslant q} \overline\chi(n) \stretchint{7ex}_{\bs I/c^6}g_{I/c^6}(n;t)t^{-\frac{3}{2}}dt+O(\log q).
\end{align*}
Since $a,b\geqslant x$ we have that
\begin{align*}
    c^3\left(\dfrac{1}{\sqrt{b}}-\dfrac{1}{\sqrt{a}}\right)\ll \sqrt{z}\cdot\frac{|a-b|}{a\sqrt{b}}\ll \sqrt{\dfrac{z}{x}}\ll 1,
\end{align*}
for which this term can be absorbed into $g_{I/c^6}(n)$ without affecting its majorant $h_{I/c^6}(n)$. In a similar fashion, we obtain
$$
\stretchint{7ex}_{\bs I/c^6} g_{I/c^6}(n;t)t^{-\frac{3}{2}} dt\ll h_{I/c^6}(n)\stretchint{7ex}_{\bs I/c^6} t^{-\frac{3}{2}}dt\ll h_{I/c^6}(n).
$$
In the case of $n^{-\frac{2}{3}}$ in place of $n^{-\frac{1}{2}}$, we carry it out in a similar way. Thus, we can write $\widetilde{g}_{I/c^6}(n)$ by combining these extra factors into $g_{I/c^6}(n)$, so that we still have $\widetilde{g}_{I/c^6}(n)\ll h_{I/c^6}(n)$ and
\begin{align*}
    \sum_{n\in I/c^6} \chi(n)f(n)=\dfrac{\tau(\chi)}{2\pi i}\sum_{|n|\leqslant q}\overline{\chi}(n)\widetilde{g}_{I/c^6}(n)+O(\log q).
\end{align*}
Therefore, we obtain the left-hand side of the inequality in Proposition \ref{prop:2.2.1} to be
\begin{align*}
    \dfrac{1}{\varphi(q)} &\sum_{\substack{\chi(q)\\ \chi\not=\chi_0}}\left| \sum_{c^6\leqslant z}\mu(c)\chi(c)^6\sum_{n\in I/c^6}\chi(n)f(n)\right|^2\\
    &= \dfrac{1}{\varphi(q)}\sum_{\substack{\chi(q)\\ \chi\not=\chi_0}}\left| \sum_{c^6\leqslant z} \mu(c)\chi(c)^6\cdot \dfrac{\tau(\chi)}{2\pi i}\sum_{|n|\leqslant q} \overline\chi(n)\widetilde{g}_{I/c^6}(n)\right|^2+O(z^{\frac{1}{3}}\log^2 q).
\end{align*}
Let $V\geqslant 0$ be a smooth function supported on $[\frac{1}{2},4]$ so that $V$ is a majorant of the characteristic function of the interval $[1,2]$. Then, we add the term corresponding to the principal character $\chi_0$, so that we can use \eqref{eq:identoties-q-aspect2} without the main term otherwise. We need to add an error term $O(z^{\frac{1}{3}})$ doing so, which is absorbed into the above error term.\\

Now, we consider double dyadic intervals in $n$ and $c$ with the majorant principle by using Equation \eqref{eq:identoties-q-aspect2}, setting $\alpha=1$, back and forth. This can be done using triangle inequality and substituting $b_n\mapsto |b_n|$ on the way back as follows. 
\begin{align*}
    \frac{1}{\varphi(q)}\sum_{\chi(q)}\left| \sum_{n}b_n\chi(n)\right|^2 &=\sum_{\ell\in(\Z/q\Z)^\times} \left|\sum_{n\equiv\ell(q)}b_n\right|^2\\
    &\leqslant \sum_{\ell\in(\Z/q\Z)^\times} \left| \sum_{n\equiv\ell(q)} |b_n|\right|^2= \frac{1}{\varphi(q)}\sum_{\chi(q)}\left| \sum_{n}|b_n|\chi(n)\right|^2.
\end{align*}
Thus, we are left to determine
\begin{align*}
\sup_{\substack{C\leqslant z^{\frac{1}{6}}\\ N\leqslant q}} &\dfrac{1}{\varphi(q)}\sum_{\chi(q)}\left|\sum_{C<c\leqslant 2C}\mu(c)\chi(c)^6\sum_{N<n\leqslant 2N} \overline\chi(n)\widetilde{g}_{I/c^6}(n)\right|^2\\
&\ll \sup_{\substack{C\leqslant z^{\frac{1}{6}}\\ N\leqslant q}} \dfrac{h^2_{I/C^6}(N)}{\varphi(q)}\sum_{\chi(q)}\left| \sum_{c}\chi(c)^6V\left(\dfrac{c}{C}\right)\sum_{n}\chi(n)V\left(\dfrac{n}{N}\right)\right|^2.
\end{align*}
Then, we split up the summation over $c$ and over $n$ by employing the Cauchy-Schwarz inequality. We can show using the decay of $\widetilde{V}$, the Mellin transform of $V$ that
$$\sum_{\chi^6\not=\chi_0}\left|\sum_{c}\chi^6(c)V\left(\dfrac{c}{C}\right)\right|^4\ll qC^2\log^4 x\text{ and }\sum_{\chi\not=\chi_0}\left|\sum_{n}\chi(n)V\left(\dfrac{n}{N}\right)\right|^4\ll qN^2\log^4 x.$$
We shall give the proof of the first inequality since the second can be done similarly. First off, we write using the similar idea to the discussion around Equation \eqref{eq:Mellin's discrete convolution} in the previous section, or we can simply use Mellin discrete convolution to write
\begin{align*}
    \sum_{c}\chi^6(c)V\left(\frac{c}{C}\right)=\frac{1}{2\pi i}\stretchint{7ex}_{\bs (\frac{1}{2})}  L(s,\chi^6)\widetilde{V}(s)C^s ds,
\end{align*}
since $\chi^6\not=\chi_0$ there is no pole at $s=1$. We mean by $\widetilde{V}(s)\ll_A (1+|\Im(s)|)^{-A}$ the decay of $\widetilde{V}(s)$, which can be proved using partial integration repeatedly. Thus, using H\"older's inequality
\begin{align*}
    \sum_{\chi^6\not=\chi_0} \left|\sum_{c}\chi^6(c)V\left(\dfrac{c}{C}\right)\right|^4 &\ll \left|\stretchint{7ex}_{\bs (\frac{1}{2})}  L(s,\chi^6)\widetilde{V}(s)C^s ds\right|^4
    \\ &\ll C^2\sum_{T\geqslant 1}{\vphantom{\sum}}' \stretchint{7ex}_{\bs T<|t|\leqslant 2T}\left|  L(\tfrac{1}{2}+it,\chi^6)\right|^4|t|^{-3} dt\ll qC^2\log^4 x,
\end{align*}
where the prime in the sum indicates the sum on the dyadic range of $T\geqslant 1$, and we have used Lemma \ref{lem:hyb-fourth-moment-subconvexity-sigma>1/2} on the fourth moment estimate. Using the Cauchy-Schwarz inequality on the above two estimates we have just proved, the bound we are interested in is now reduced to
\begin{align*}
&\ll (\log^6 x) \sup_{\substack{C\leqslant z^{1/6}\\ N\leqslant q}} h_{I/C^6}(N)^2qCN\\
&\ll (\log^6 x)\cdot\left(\sup_{\substack{C\leqslant \sqrt[6]{|I|/q}}}qC+\sup_{\substack{\sqrt[6]{|I|/q}\leqslant C\leqslant z^{1/6}\\}}\sqrt[6]{|I|/q}\right)\ll |I|^{\frac{1}{6}}q^{\frac{5}{6}}(\log^6 x).
\end{align*}
This completes our proof of Proposition \ref{prop:2.2.1}.
\end{proof}

Let $q$ be a prime. From now on, we may introduce the set of quadratic nonresidue modulo $q$, 
\begin{align}
    \label{eq:QNR-def}
    \text{QNR}(q):=\{\alpha\in(\Z/q\Z)^\times: \alpha\text{ is a quadratic nonresidue modulo }q\}.
\end{align}

\begin{proposition}
\label{prop:2.2.2}
Let $\eps\in(0,\frac{1}{1000})$ and $x^{\frac{1}{11}}\leqslant q\leqslant x^{1-50\eps}$ be a prime. Let $g(u)$ be a smooth function with support on $[1,2]$ that is $1$ for $u\in [1+(x/q)^{-\eps/4}, 2(1-(x/q)^{-\eps/4})]$, and for $k\geqslant 0$, we have $g^{(k)}(u)\ll (x/q)^{k\eps/4}$. Let $(x/q)^{1+\eps}\leqslant z\leqslant x^{-\eps}\sqrt{qx}$. Then, for any $\alpha\in\textup{QNR}(q)$, we obtain
    \begin{align}
    \label{eq:lhs-prop-3.4}
        \dfrac{1}{\varphi(q)}\sum_{\substack{\chi(q)\\\chi\not=\chi_0}}\dfrac{1+\chi(\alpha)}{2}\left|\sum_{\substack{n\geqslant 1\\ c^6\leqslant z}}\mu(c)\chi(n)\chi(c^6)f(n)g\left(\frac{n}{x/c^6}\right)\right|^2\sim Cx^{\frac{1}{6}}q^{\frac{5}{6}}.
    \end{align}
\end{proposition}

\begin{proof}
    Let $h(n(x/c^6)^{-1})=f(n)g(n(x/c^6)^{-1})$. We use the Poisson summation for characters,
    $$\sum_{n} \chi(n)h\left(\dfrac{n}{x/c^6}\right)=\tau(\chi)\cdot\dfrac{x}{qc^6}\sum_{\ell}\overline\chi(\ell)\widehat{h}(\dfrac{x\ell}{qc^6}).$$
    Thus, the left-hand side of \eqref{eq:lhs-prop-3.4} is 
    \begin{align*}
        &\dfrac{1}{\varphi(q)} \sum_{\chi\not=\chi_0}\frac{1+\chi(\alpha)}{2}\left| \sum_{c^6<z}\mu(c)\chi^6(c)\sum_{n}\chi(n){h}(\dfrac{n}{x/c^6})\right|^2
        \\ &= \sum_{\chi\not=\chi_0}\frac{1+\chi(\alpha)}{2\varphi(q)}\left| \sum_{c^6<z}\mu(c)\chi^6(c)\tau(\chi)\cdot \dfrac{x}{qc^6}\sum_{\ell}\overline\chi(\ell) \widehat{h}(\dfrac{x\ell}{qc^6})\right|^2
        \\
        &= \frac{q}{\varphi(q)}\dfrac{x^2}{q^2}\sum_{\chi\not
        =\chi_0}\frac{1+\chi(\alpha)}{2}\left|\sum_{c^6\leqslant z, \ell\in \mathbb Z}\frac{\mu(c)}{c^6}\chi(c^6)\overline \chi(\ell) \widehat h(\frac{x\ell}{c^6q})\right|^2
    \end{align*}
which is
\begin{align}
    \label{term:arith-after-expanding}
    \frac{x^2}{2q}\sum_{\substack{n_1,n_2\in\mathbb Z\\c_1^6,c_2^6\leqslant z\\n_1c_1^6\equiv n_2c_2^6,\alpha n_2c_2^6(q)\\\gcd(q,cc_1n_1n_2)=1}}\frac{\mu(c_1)}{c_1^6}\frac{\mu(c_2)}{c_2^6} \widehat{h}(\frac{xn_2}{c_1^6q})\overline{ \widehat{h}(\frac{xn_1}{c_2^6 q})}+O(z^{1/3}x^{2\varepsilon/3}),
\end{align}
where the error term arises from the summand $\chi=\chi_0$ as follows. We know that
\begin{align*}
    \widehat{h}(\xi) =-\stretchint{7ex}_{\bs \R} \frac{e^{2\pi iu\xi}}{2\pi i\xi}f'(u) du \ll |\xi|^{-1}(x/q)^{\frac{\eps}{4}},
\end{align*}
and whence,
\begin{align*}
    \frac{x^2}{q}\sum_{\substack{n_1,n_2\in\mathbb Z\\c_1^6,c_2^6\leqslant z\\n_1c_1^6\equiv n_2c_2^6,\alpha n_2c_2^6(q)\\\gcd(q,cc_1n_1n_2)=1}} \frac{1}{c_1^6c_2^6} \left| \widehat{h}(\frac{xn_2}{c_1^6q})\overline{ \widehat{h}(\frac{xn_1}{c_2^6 q})}\right| \ll \frac{x^{2+\eps/2}}{q} \sum_{\substack{n_1,n_2\in\mathbb Z^+\\c_1^6,c_2^6\leqslant z\\n_1c_1^6\equiv n_2c_2^6,\alpha n_2c_2^6(q)}} \frac{1}{c_1^6c_2^6} \cdot\frac{c_1^6c_2^6q^2}{x^2n_1n_2}.
\end{align*}
We can add or drop conditions $|n_1|,|n_2|\leqslant x^{\varepsilon/3}zq/x$ and $c,c_1\geqslant x^{1/6-\varepsilon/18}/q^{1/6}$ as we want, since the error terms occurring are acceptable. This can be done by the same trick of repeatedly partial integrations on $e(\cdot)$ in $ \widehat{h}$. Separate the set of all possible $(n_1,n_2)$ into $2$ possibilities, those $(n_1,n_2)$ that are in
$$
\mathcal M:=\{(k_1^2m,k_2^2m):m\text{ is square-free},k_1,k_2\in\mathbb N\},
$$
and those that are not in this set.

\subsubsection{The main term $(n_1,n_2)\in\mathcal{M}$}
We have that $mk_1^2c_1^6\equiv mk_2^2c_2^6 \pmod{q}$ and $\gcd(q,n_1n_2)=1$ implies $k_1^2c_1^6\equiv k_2^2c_2^6\pmod{q}$. The case $mk_1^2c_1^6\equiv\alpha mk_2^2c_2^6\pmod{q}$ implies no solution since $\alpha\in\textup{QNR}(q)$, so in this case it does not produce a main term. Thus, for the former case, $k_1c_1^3=k_2c_2^3$ since $k_1c_1^3+k_2c_2^3\ll \sqrt{x^{\eps/3}zq/x}\sqrt{z}\leqslant o(q)$ since $z<x^{-\varepsilon}\sqrt{qx}$. By this condition, we parametrize

$$k_1=\dfrac{c_2^3\ell}{(c_1,c_2)^3}\text{ and }k_2=\frac{c_1^3\ell}{(c_1,c_2)^3} \text{ : }(\ell\geqslant 1).$$

The contribution from such $n,n_1$, by replacing $m\ell^2$ with simply $\ell$ since any positive integer $n$ can be written uniquely as $m\ell^2$ for square-free $m$, is 

$$
\dfrac{x^2}{q}\sum_{c_1^6,c_2^6\leqslant z, (q,cc_1)=1}\dfrac{\mu(c_1)\mu(c_2)}{c_1^6c_2^6}\sum_{\ell\geqslant 1, (\ell,q)=1}\left| \widehat{h}(\frac{x\ell}{q(c_1,c_2)^6})\right|^2.
$$

Let $H=x/q$ and write $W= \widehat h$. Consider from the condition on $f,g$ by bounding $f\ll 1$ that
$$
\sum_{\lambda\geqslant 1}|W(\lambda/\nu)|^2\ll \sum_{\lambda\leqslant \nu H^{\varepsilon/4}}1+\sum_{\lambda>\nu  H^{\varepsilon/4}} \dfrac{H^{\varepsilon/2}\nu^2}{\lambda^2}\ll \nu H^{\varepsilon/4}.
$$
Therefore, by letting $\lambda=\ell, \nu=H^{-1}(c_1,c_2)^6$,
\begin{align*}
    H^2&\sum_{\substack{c_1>z^{1/6}||\\ c_2>z^{1/6}}} \dfrac{\mu(c_1)\mu(c_2)}{c_1^6c_2^6}\sum_{\ell\geqslant 1}\left| \widehat{h}(\frac{x\ell}{q(c_1,c_2)^6})\right|^2 \ll H^{1+\varepsilon/4} \sum_{\substack{c_1>z^{1/6}||\\ c_2>z^{1/6}}}\dfrac{(c_1,c_2)^6}{c_1^6c_2^6}
    \\ &\ll H^{1+\varepsilon/4}\sum_{c_0\geqslant 1}\frac{1}{c_0^6}\sum_{\gamma\geqslant H^{(1+2\varepsilon)/6}/c_0,\gamma_1\geqslant 1} \dfrac{1}{\gamma^6\gamma_1^6}\ll H^{1+\varepsilon/4}\sum_{c_0\geqslant 1} \dfrac{1}{c_0^6}\min(1,\frac{c^5_0}{H^{5(1+2\varepsilon)/6}})\ll H^{1/6-\varepsilon/20}
\end{align*}
when $z\geqslant H^{1+2\varepsilon}$, which satisfies our condition on $z$. Here, $P||Q$ means the statement holds if $P$ or $Q$ holds. Since $q\geqslant x^{\frac{1}{11}}$ and $c_1,c_2\leqslant z^{\frac{1}{6}}\leqslant x^{-\eps/6}(qx)^{\frac{1}{12}}\leqslant o(q)$, we may drop the condition $\gcd(q,c_1c_2)=1$. Similarly, for $\ell\geqslant q$, the contribution is $O_{A}(x^{-A})$ by the inequality $\widehat{h}(\xi)\ll_{A,\eps} x^{-A}$ so it can be negligible. We are left to study 

$$H^2\sum_{\ell\geqslant 1}\sum_{c_1,c_2\geqslant 1}\dfrac{\mu(c_1)\mu(c_2)}{c_1^6c_2^6}\left|W(\frac{H\ell}{(c_1,c_2)^6})\right|^2+O(H^{1/6-\varepsilon/20}).$$

Let $g(x):=|W(e^x)|^2e^x$ we have that for $r>0$ and $\Re s=c\in(-1,1)$,
$$
|W(r)|^2=\frac{r^{-1}}{2\pi i}\stretchint{7ex}_{\bs(c)}r^s \widehat g(\frac{s}{2\pi i})ds.
$$
Then, the main term above becomes
\begin{align*}
&\dfrac{H}{2\pi i} \sum_\ell\sum_{c_1,c_2}\dfrac{\mu(c_1)\mu(c_2)\gcd(c_1,c_2)^6}{c_1^6c_2^6}\dfrac{1}{\ell}\stretchint{7ex}_{\bs(c)}H^s\ell^{s}\gcd(c_1,c_2)^{-6s} \widehat g(\frac{s}{2\pi i})ds
\\ &\quad =\dfrac{H}{2\pi i} \stretchint{7ex}_{\bs(c)} H^s\zeta(1-s)\sum_{c_1,c_2\geqslant 1} \dfrac{\mu(c_1)\mu(c_2)\gcd(c_1,c_2)^{6-6s}}{c_1^6c_2^6} \widehat g(\frac{s}{2\pi i})ds.
\end{align*}
Now we compute the sum that occurred in the integrand. 
\begin{align*}
\sum_{c_1,c_2\geqslant 1} &\dfrac{\mu(c_1)\mu(c_2)}{c_1^6c_2^6}\gcd(c_1,c_2)^{6-6s}\\
&= \sum_{c\geqslant 1} \dfrac{\mu^2(c)}{c^{6+6s}}\sum_{(\gamma_1,c)=1}\dfrac{\mu(\gamma_1)}{\gamma_1^6}\sum_{\substack{(\gamma_2,c)=1\\ (\gamma_2,\gamma_1)=1}} \dfrac{\mu(\gamma_2)}{\gamma_2^6}\\
&= \prod_{p}\left(1-\dfrac{2}{p^6}\right)\sum_{c\geqslant 1}\dfrac{\mu^2(c)}{c^{6+6s}}\prod_{p|c}\left(1-\dfrac{2}{p^6}\right)^{-1}=\prod_{p} \left(1-\dfrac{2}{p^6}+\dfrac{1}{p^{6+6s}}\right).
\end{align*}
Plugging in the calculation and splitting the term $\zeta(6+6s)$ from it, the above main term is now
\begin{align*}
\frac{H}{2\pi i}&\stretchint{7ex}_{\bs (c)} H^s\zeta(1-s)\zeta(6+6s)\prod_{p}\left(1-\dfrac{2}{p^6}+\dfrac{2}{p^{12+6s}}-\dfrac{1}{p^{12+12s}}\right)  \widehat{g}(\frac{s}{2\pi i})ds\\
&= \frac{H^{\frac{1}{6}}}{2}\zeta(\frac{11}{6})\prod_p\left(1-\dfrac{1}{p^2}+\dfrac{2}{p^{6}}+\dfrac{2}{p^7}\right)\stretchint{7ex}_{\bs 0}^{\infty}|W(y)|^2y^{\frac{5}{6}}dy +O(H^{\frac{1}{6}-\frac{\eps}{20}}),
\end{align*}
by moving the contour to $\Re s=-\frac{11}{12}$. Note that the singularity at $s=-\frac{5}{6}$ is our desired main term, and we are done for these diagonal terms $(n_1,n_2)\in\mathcal{M}$.

\subsubsection{The term $(n_1,n_2)\not\in\mathcal{M}$}

As we mentioned, we might impose the condition $c_1,c_2\geqslant x^{1/6-\eps/20}/q^{1/6}$ as desired. Using dyadic ranges over $C_1,C_2\in [x^{1/6-\eps/20}/q^{1/6}, z^{1/6}]$ we need to study
\begin{align}
\label{eq:arith-casenotin-M}
    \sup_{\ell\in(\Z/q\Z)^\times}\dfrac{x^2}{q}\cdot\dfrac{1}{C_1^6C_2^6}\sum_{(n_1,n_2)\not\in\mathcal{M}} V(\frac{n_1}{N_1})V(\frac{n_2}{N_2})\sum_{\substack{n_1c_1^6\equiv \ell n_2c_2^6(q)\\ (q,n_1n_2c_1c_2)=1}} V(\frac{c_1}{C_1})V(\frac{c_2}{C_2}),
\end{align}
for some smooth function $V$ (absorbing $ \widehat{h}$) whose support is in $[1/4,4]$. Using characters, now the Term \eqref{eq:arith-casenotin-M} is bounded (without using the condition not in $\mathcal{M}$) by the constant multiple of
\begin{align*}
    \dfrac{x^2}{q\varphi(q)}&\cdot\dfrac{1}{C_1^6C_2^6}\sum_{\chi^6\not=\chi_0} \left(\sum_{n_1}\chi(n_1)V(\frac{n_1}{N_1})\right)\left(\sum_{n_2}\overline\chi(n_2)V(\frac{n_2}{N_2})\right)\\
    &\times \left(\sum_{c_1}\chi^6(c_1)V(\frac{c_1}{C_1})\right)\left(\sum_{c_2}\overline\chi^6(c_2)V(\frac{c_2}{C_2})\right)+O\left(\dfrac{x^2}{q^2}\dfrac{N_1N_2}{C_1^5C_2^5}\right).
\end{align*}
where the error term is negligible since $N_j\ll x^{\eps/3}C_j^6q/x$ for $j=1,2$. By the contour integral and H\"older inequality, the main term is bounded by
\begin{align*}
    \dfrac{x^2}{q^2}\cdot\dfrac{\sqrt{N_1N_2C_1C_2}}{C_1^6C_2^6}\sum_{\chi}\stretchint{7ex}_{\bs |u|\leqslant x^{\eps/3}} |L(\tfrac{1}{2}+iu,\chi)|^4 du.
\end{align*}
When $C_1^3C^3_2\geqslant x^{1+2\eps}/q$ the above becomes $$\dfrac{x^2}{q}x^{\eps/2}\dfrac{\sqrt{N_1N_2C_1C_2}}{C_1^6C_2^6}\ll x^{5\eps/6}x\dfrac{1}{(C_1C_2)^{5/2}}\ll x^{\frac{1}{6}-\frac{\eps}{10}}q^{\frac{5}{6}}$$
by the fourth moment bound of the $ L$-function (Lemma \ref{lem:hyb-fourth-moment-subconvexity-sigma>1/2}) and the lower bound on $C_1C_2$ we assumed.\\

It is left to show a similar result when $C_1C_2\leqslant x^{1/3+2\eps/3}/q^{1/3}$ by bringing back the conditions that $c_1,c_2$ are square-free. Since $C_1,C_2\geqslant (x/q)^{1/6}x^{-\eps/20}$, we have $C_1,C_2\leqslant (x/q)^{1/6}x^{2\eps}$ and so $N_1,N_2\ll x^{3\eps}$. When $(n_1,n_2)\not\in \mathcal{M}$ and $n_1c_1^6\equiv n_2c_2^6\pmod{q}$ we have to count solutions to $n_1c_1^6=n_2c_2^6+q\ell$ with $0<|\ell|\ll x^{1+3\eps}/q^2$, since $\ell=0$ and $c_1,c_2$ are square-free, implying $n_1=n_2$ and $(n_1,n_2)\in\mathcal{M}$ with $k_1=k_2=1$. The number of such solutions when fixing $n_1,n_2,\ell$ is $O_\eps(x^{\eps})$. Totally, the number of $(n_1,n_2,c_1,c_2,\ell)$ in this case is at most $O_\eps(x^{1+9\eps}/q^2)$. Note that we need to use the bound on $ \widehat{h}(u)\ll (ux/z)^{-\frac{1}{2}}$ to acquire the last term on the left-hand side below by considering the Term \eqref{term:arith-after-expanding}. The contribution in this case is therefore bounded by
\begin{align}
    \dfrac{x^2}{q}\cdot\dfrac{1}{C_1^6C_2^6}\cdot\dfrac{x^{1+9\eps}}{q^2}\times \left(\dfrac{\sqrt
    qz}{x}\right)^2\ll qx^{8\eps}\ll x^{1/6}q^{5/6}
\end{align}
where the last inequality is true for $q\ll x^{1-50\eps}$ as desired. The condition $n_1c_1^6\equiv \alpha n_2c_2^6\pmod{q}$ can be done similarly with the implied constant possibly dependent on the fixed $\alpha$.
\end{proof}

Now, we combine Propositions \ref{prop:2.2.1} and \ref{prop:2.2.2} to prove Proposition \ref{prop:equiv-to-gcd-smallc}. As we noted, we shall fill in the gap in a smooth function $g$ (recall that $\text{supp}(g)\subset [1,2]$) using appropriate bounds in small intervals $I_1/x:=[1, 1+(x/q)^{-\eps/4}], I_2/x:=[2(1-(x/q)^{-\eps/4}),2]$, since we have $g\equiv 1$ in between. This is when we use Proposition \ref{prop:2.2.1} since the bound depends on the size of the intervals, which are small in this case.\\

It suffices to show that
\begin{align}
    \label{eq:suffice to show in Prop 3.1}
    \frac{1}{\varphi(q)}\sum_{\chi\not=\chi_0}\left|\sum_{\substack{c^6\leqslant z\\ \frac{n}{x/c^6}\in \frac{I_1}{x}\cup \frac{I_2}{x}}}\mu(c)\chi(n)\chi(c^6)f(n)\left(1-g(\frac{n}{x/c^6})\right)\right|^2 \ll x^{\frac{1}{6}}q^{\frac{5}{6}}x^{-\frac{\eps}{16}},
\end{align}
and we will prove on the condition $n(x/c^6)^{-1}\in I_1/x$, since the other case proof is similar. We shall be considering
\begin{align*}
   \sum_{\substack{m\in I_1}} \chi(m)\left(1-g(\frac{m}{x})\right)\sum_{\substack{m=nc^6\\ c^6\leqslant z}} \mu(c)f(n)=:\sum_{\substack{m\in I_1}} \chi(m)\left(1-g(\frac{m}{x})\right)M_m.
\end{align*}
Here, we write the inner sum above as $M_m:=M_m^{(z)}$. Proposition \ref{prop:2.2.1} corresponds to the bound involving the above except the term $(1-g(m/x))$. Now, we employ partial summation (to split the term $1-g(n/x)$ out so that we can use Proposition \ref{prop:2.2.1}) and the Cauchy-Schwarz inequality to obtain that
\begin{align*}
    \frac{1}{\varphi(q)}&\sum_{\chi\not=\chi_0}\left|\sum_{n\in I_1} \chi(n)M_n\cdot \left(1-g(\frac{n}{x})\right)\right|^2
    \\ &\ll \frac{1}{\varphi(q)}\sum_{\chi\not=\chi_0}\left|\sum_{n\in I_1}\chi(n)M_n\right|^2+\frac{(x/q)^{\frac{\eps}{4}}}{x}\stretchint{7ex}_{\bs [x,x+x(x/q)^{-\eps/4}]} \frac{1}{\varphi(q)}\sum_{\chi\not=\chi_0}\left| \sum_{n\leqslant u}\chi(n)M_n\right|^2 du.
\end{align*}
Invoking Proposition \ref{prop:2.2.1}, we obtain Equation \eqref{eq:suffice to show in Prop 3.1} since $|I_1|, u\ll x(x/q)^{-\eps/4}$, and hence, the proof of Proposition \ref{prop:equiv-to-gcd-smallc} is now completed.

\subsection{Proof of Proposition \ref{prop:equiv-to-gcd-largec}}
Firstly, we bound the left-hand side by a constant multiple of
\begin{align*}
    \dfrac{1}{\varphi(q)^2} &\sum_{\chi\not=\chi_0}\left| \sum_{\substack{nc^6\leqslant x\\ C<c\leqslant 2C}}\mu(c)\chi(n)\chi^6(c)f(n)\right|^2\\
    &\ll x^{\frac{1}{3}}\log^{10}x \sup_{\substack{z\leqslant C^6\leqslant x\\ T\leqslant x^3}} \dfrac{1}{\varphi(q)^2}\sum_{\chi\not=\chi_0}\dfrac{1}{T}\stretchint{7ex}_{\bs |t|\leqslant T} | L(\tfrac{2}{3}+it,\chi)|^2| \mathcal M(\tfrac{1}{6}+it,\chi^6)|^2 dt,
\end{align*}
using the Perron formula and moving the path of the contour to $\Re s=1/6$ with the acceptable error term from the subconvexity bound, where
$$\mathcal M(\tfrac{1}{6}+it,\chi^6):=\sum_{C<c\leqslant 2C} \dfrac{\mu(c)\chi^6(c)}{c^{1+6it}}\ll 1.$$
For $V>0$, let us write
$$S_{T,q}(V):=\{(t,\chi): V\leqslant |\mathcal M(\tfrac{1}{6}+it,\chi^6)|\leqslant 2V, |t|\leqslant T, \chi\,(\text{mod } q)\}.$$
It suffices to show that

\begin{equation}
    \label{eq:bdd-on-L-func}
    \dfrac{x^{\frac{1}{3}}V^2}{\varphi(q)}\sum_{\chi\not=\chi_0} \dfrac{1}{T}\stretchint{7ex}_{\bs t: (t,\chi)\in S_{T,q}(V)} | L(\tfrac{2}{3}+it, \chi)|^2dt \ll x^{\frac{1}{6}}q^{\frac{5}{6}}(x/q)^{-\eps/10}.
\end{equation}

Split the interval $[-T,T]$ into intervals $I$ of length $1$, then if $I\cap S_{T,q}(V)$ is nonempty, we choose a point in $I\cap S_{T,q}(V)$. Write $\widetilde{S}_{T,q}(V)$ as a collection of the points in question for the odd-positioned intervals so that $\widetilde{S}_{T,q}(V)$ is well-spaced. Thus, we have captured the Lebesgue measure of $S_{T,q}(V)$ within $\widetilde{S}_{T,q}(V)$, so
$$\lambda(S_{T,q}(V))\ll \big|\widetilde{S}_{T,q}(V)\big|\ll V^{-2}+qT\min(C^{-1}V^{-2}, C^{-2}V^{-6}),$$
in which the last inequality is justified by Lemma \ref{lem:Hybrid large-value theorem}. When the first term dominates, we have that the left-hand term in Equation \eqref{eq:bdd-on-L-func} is
\begin{align*}
    \ll \dfrac{x^{1/3+\eps}}{qT}(qT)^{1/3+\eps}\ll \dfrac{x^{1/3+\eps}}{q^{2/3-\eps}}\ll x^{1/6}q^{5/6}(x/q)^{-\eps/10}
\end{align*}
when $q\gg x^{1/9+\eps}$. Now, we consider when the second term dominates. In this case, we employ the Cauchy-Schwarz inequality and the inequality $\min(a,b)\leqslant \sqrt{ab}$ to obtain the bound to be
\begin{align*}
    &\dfrac{x^{\frac{1}{3}}V^2}{qT} \lambda(S_{T,q}(V))^{\frac{1}{2}} \left(\sum_{\chi\not=\chi_0} \stretchint{7ex}_{\bs |t|\leqslant T} | L(\tfrac{2}{3}+it)|^4dt\right)^\frac{1}{2}\\
    &\ll x^{1/3}(\log^2) xV^2\min(C^{-1}V^{-2}, C^{-2}V^{-6})^{1/2}\leqslant x^{1/3}(\log^2 x)C^{-3/4}\ll x^{1/6}q^{5/6}(x/q)^{-\eps/10} 
\end{align*}
since $C^6\gg z\geqslant (x/q)^{1+\eps}$ and $q\gg x^{\frac{1}{23}+\eps}$. This finishes the proof of Proposition \ref{prop:equiv-to-gcd-largec}.\\

Let us combine Propositions \ref{prop:equiv-to-gcd-smallc} and \ref{prop:equiv-to-gcd-largec} together. This implies, using the Cauchy-Schwarz inequality, that
\begin{equation}
\label{eq:eq-combine-4.1,4.2}
    \dfrac{1}{\varphi(q)}\sum_{\ell\in(\Z/q\Z)^\times}\left|\frac{1}{2}\sum_{\substack{x<nc^6\leqslant 2x\\ nc^6\equiv\ell,\alpha\ell(q)}}\mu(c)f(n)-\dfrac{1}{\varphi(q)}\sum_{\substack{x<nc^6\leqslant 2x\\ \gcd(nc^6,q)=1}}\mu(c)f(n)\right|^2\sim C(x/q)^{\frac{1}{6}},
\end{equation}
for all $x^{\frac{1}{9}+\eps}\leqslant q\leqslant x^{1-50\eps}$. To summarize, Equation \eqref{eq:eq-combine-4.1,4.2} approximates the sum over the condition $nc^6\equiv\ell\pmod{q}$ by the sum over the condition $\gcd(nc^6,q)=1$, which is more comfortable to produce a main term counting square-full integers in arithmetic progressions.\\

Now, we state the proposition that handles the residual term $r(n)$. Proposition \ref{prop:err-term} below is done using the bound on the $(2,3)$-divisor problem, i.e., Equation \eqref{eq:err-23-div} and similar counting as in the proof of Proposition \ref{prop:equiv-to-gcd-smallc}.

\begin{proposition}
\label{prop:err-term}
Let $\eps\in(0,\frac{1}{1000})$ be given. For a prime $x^{\frac{51}{114}+\eps}\leqslant q\leqslant x^{1-\eps}$ with $\alpha\in\textup{QNR}(q)$ be any quadratic nonresidue modulo $q$, we have that
\begin{align*}
    \dfrac{1}{\varphi(q)}\sum_{\ell\in (\mathbb Z/q\mathbb Z)^\times}\left|\sum_{\substack{x<nc^6\leqslant 2x\\ nc^6\equiv\ell,\alpha\ell (q)}}\mu(c)r(n)\right|^2\ll_\eps (x/q)^{\frac{1}{6}-\frac{\eps}{16}}.
\end{align*}
\end{proposition}

Also, we need to ensure that the term
\begin{align*}
    \dfrac{1}{\varphi(q)}&\sum_{\substack{nc^6\leqslant x\\ (q,nc^6)=1}}\mu(c)r(n)=\dfrac{1}{\varphi(q)}\sum_{\substack{c^6\leqslant x\\ (c,q)=1}}\mu(c)\sum_{\substack{n\leqslant x/c^6\\ (n,q)=1}}r(n),
\end{align*}
is also acceptable. The inner sum above is
\begin{align*}
    \sum_{n\leqslant x/c^6} &r(n)-\sum_{\substack{n\leqslant x/c^6\\ q|n}} r(n)\ll \dfrac{x^{\frac{2}{15}}}{c^{\frac{4}{5}}}+\dfrac{x^{\frac{1}{2}}}{qc^3},
\end{align*}
where we invoke Equation \eqref{eq:err-23-div} for the first sum, and we consider whether $n$ is of the form $a^2b^3$ or not for the second. Plugging in this bound and summing over $c^6\leqslant x$, we obtain that its contribution is $\ll (x/q)^{1/6-\eps/10}$ for $q\gg x^{3/11+5\eps}$. Combining Proposition \ref{prop:err-term} with the above calculation, we obtain that
\begin{equation}
\label{eq:eq-4.5}
    \dfrac{1}{\varphi(q)}\sum_{\ell\in(\Z/q\Z)^\times} \left|\frac{1}{2}\sum_{\substack{x<nc^6\leqslant 2x\\ nc^6\equiv\ell,\alpha\ell(q)}}\mu(c)r(n)-\dfrac{1}{\varphi(q)}\sum_{\substack{x<nc^6\leqslant 2x\\ \gcd(nc^6,q)=1}}\mu(c)r(n)\right|^2\ll (x/q)^{\frac{1}{6}-\frac{\eps}{16}},
\end{equation}
for all $q\geqslant x^{51/114+\eps}$.

\begin{proof}[Proof of Proposition \ref{prop:err-term}]

Recall that $r(n)=1_{2,3}(n)-\frac{\zeta(\frac{3}{2})}{2{n}^{\frac{1}{2}}}-\frac{\zeta(\frac{2}{3})}{3n^{\frac{2}{3}}}$. We obtain by the orthogonality that
\begin{align}
\label{term:arith-var-of-err}
    &\sum_{\ell\in(\mathbb Z/q\mathbb Z)^\times} \left|\sum_{\substack{nc^6\leqslant x\\nc^6\equiv \ell,\alpha\ell(q)}}\mu(c)r(n)\right|^2
    \\ &\qquad=\sum_{\substack{c_1^6,c_2^6\leqslant x\\\gcd(q,c_1c_2n_1n_2)=1}} \mu(c_1)\mu(c_2)\sum_{n_1\leqslant x/c_1^6}r(n_1)\sum_{\substack{n_2\leqslant x/c_2^6\\n_2\equiv n_1c_1^6c_2^{-6},\alpha n_1c_1^6c_2^{-6}(q)}}r(n_2). \nonumber
\end{align}
It suffices to consider for fixed $c_1,c_2\leqslant x^{1/6}$,

$$\sum_{n_1=a^2b^3\leqslant x/c_1^6}\sum_{\substack{n_2\leqslant x/c_2^6\\ n_2\equiv n_1c_1^6c_2^{-6}(q)||\\  n_2 \equiv \alpha n_1c_1^6c_2^{-6}(q)}}r(n_2)\text{ and }\sum_{n_1\leqslant x/c_1^6}h(n)\sum_{\substack{n_2\leqslant x/c_2^6\\ n_2\equiv n_1c_1^6c_2^{-6}(q)||\\  n_2 \equiv \alpha n_1c_1^6c_2^{-6}(q)}}r(n_2),$$

\noindent
by splitting $r(n_1)$, where $h(n)= n^{-1/2}$ or $n^{-2/3}$. The first double sum needs the use of quadratic nonresiduality of $\alpha$, and we shall treat it later. The second double sum can be bounded dyadically on $n_1$ without the use of $\alpha$, i.e., for fixed $y\in [1, x/c_1^6]$, the contribution is at most

\begin{align}
    \label{eq:second-dobule-sum-err-prop}
    \sum_{c_1,c_2\leqslant x^{1/6}}\left|\sum_{y<n_1\leqslant 2y}{\vphantom{\sum}}' \dfrac{1}{n_1^{\frac{1}{2}}}\sum_{\substack{n_2\leqslant x/c_2^6\\ n_2\equiv n_1c_1^6c_2^{-6}(q)}}r(n_2)\right|
\end{align}

\noindent
where the prime denotes the sum over $n_1$ in that range with a multiple of $q$ terms of $n_1$. We replace $h(n)=n^{-1/2}$ now, since the other case can be done similarly with slightly easier arguments. Let us write $\delta=\tfrac{8}{255}$. In this case, we shall split into the case when $y\in [q^{\frac{7}{5}+\delta}, x/c_1^6]$ and $y\in [1, q^{\frac{7}{5}+\delta}]$. The residual term, i.e., the terms in \eqref{eq:second-dobule-sum-err-prop} in $[y,2y]$ outside the primed sum, is bounded by $\log x$ times
$$\sum_{c_1,c_2\leqslant x^{1/6}}\sum_{n_1\leqslant q}\dfrac{1}{n_1^{\frac{1}{2}}}\dfrac{x^{1/2}}{qc_2^3}\ll x^{1/6}\sqrt{\dfrac{x}{q}}\ll x^{1/6}q^{5/6}(x/q)^{-\eps},$$
whenever $q\gg x^{3/8+5\eps}$, which is acceptable. Thus, we are left to deal with the Term \eqref{eq:second-dobule-sum-err-prop} above. When $y\in [q^{\frac{7}{5}+\delta}, x/c_1^6]$, we consider the $O(y/q)$ terms of
\begin{align*}
    \sum_{a<n_1\leqslant a+q} &\frac{1}{n_1^{\frac{1}{2}}}\sum_{\substack{n_2\leqslant x/c_2^6\\ n_2\equiv n_1c_1^6c_2^{-6}(q)}}r(n_2) = \dfrac{1}{\sqrt{a+q}}\sum_{n\leqslant x/c_2^6} r(n_2)+O\left(\sum_{j\leqslant q} \left(\frac{1}{\sqrt{a+j}}-\frac{1}{\sqrt{a+q}}\right)\frac{x^{\frac{1}{2}}}{c_2^3q}\right)
\end{align*}
where $a\gg y\gg q^{\frac{7}{5}+\delta}$ is an integer, and we have used the trivial bound on the sum over $r(n_2)$. For the first term we use \eqref{eq:err-23-div}, and thus when we sum over $c_1,c_2\leqslant x^{1/6}$ the above case contributes to at most
\begin{align*}
    \sum_{c_1,c_2\leqslant x^{1/6}} &\sup_{y\in [q^{\frac{7}{5}+\delta}, x/c_1^6]}\dfrac{y}{q}\left(\frac{x^{\frac{2}{15}}}{c_2^{\frac{4}{5}}a^{\frac{1}{2}}}+\frac{x^{\frac{1}{2}}q}{c_2^3a^{\frac{3}{2}}}\right)\\
    &\leqslant\sum_{c_1,c_2\leqslant x^{1/6}} \frac{x^{\frac{2}{15}}}{qc_2^{\frac{4}{5}}}\left(\frac{x}{c_1^6}\right)^{\frac{1}{2}}+\dfrac{x^{\frac{1}{2}}}{c_2^3q^{\frac{3}{2}(\frac{7}{5}+\delta)}}\frac{x}{c_1^6}\ll \dfrac{x^{\frac{1}{2}+\frac{2}{15}+\frac{1}{30}}}{q}+\dfrac{x^{1+\frac{1}{2}}}{q^{\frac{21}{10}+\frac{3\delta}{2}}}\ll_\eps x^{1/6}q^{5/6}(x/q)^{-\eps/10},
\end{align*}
where the last inequality holds for $q\gg \max(x^{\frac{3}{11}+\eps},x^{\frac{4}{3}\frac{1}{44/15+3\delta/2}+\eps})=x^{\frac{51}{114}+\eps}$ with the second term dominates in the maximum function. When $q\in [1, q^{\frac{7}{5}+\delta}]$, we bound the Term \eqref{eq:second-dobule-sum-err-prop} trivially, i.e.,
\begin{align}
\label{eq:trivialboundon_r(n)}
    \sum_{\substack{n\leqslant y\\ n\equiv \ell (q)}}|r(n)|\ll \sum_{\substack{n=a^2b^3\leqslant y\\ n\equiv \ell(q)}} 1+\sum_{\substack{n\leqslant y\\ n\equiv \ell(q)}} \frac{1}{n^{\frac{1}{2}}}\ll\frac{y^{\frac{1}{2}}}{q},
\end{align}
uniformly on $\ell\in (\Z/q\Z)^\times$, to obtain the bound
\begin{align*}
    \sum_{c_1,c_2\leqslant x^{1/6}} \sqrt{y}\cdot \dfrac{x^{\frac{1}{2}}}{c_2^3q} \ll \dfrac{x^{\frac{1}{2}+\frac{1}{6}}\cdot q^{\frac{7}{10}+\frac{\delta}{2}}}{q}\ll x^{1/6}q^{5/6}(x/q)^{-\eps/10},
\end{align*}
for all $q\geqslant x^{\frac{1}{2}\frac{1}{17/15-\delta/2}+\eps}=x^{\frac{51}{114}+\eps}$.\\

Now, we turn to the first double sum, i.e., we determine
\begin{align*}
    \sum_{n_1=a^2b^3\leqslant x/c_1^6}\sum_{\substack{n_2\leqslant x/c_2^6\\ n_2\equiv n_1c_1^6c_2^{-6}(q)||\\  n_2 \equiv \alpha n_1c_1^6c_2^{-6}(q)}}r(n_2).
\end{align*}
We sum indexing $b$ first and split the case when $b\leqslant \frac{x^{1/3}}{10c_1^2q^{2/3}}$ or not. This is because we will obtain for the former case $\frac{x^{1/2}}{c_1^3b^{3/2}}\geqslant 10q$. For the former case, we determine
\begin{align*}
    \sum_{b\leqslant \frac{x^{1/3}}{10c_1^2q^{2/3}}}&\sum_{a\leqslant \frac{x^{1/2}}{c_1^3b^{3/2}}}{\vphantom{\sum}}'\sum_{\substack{n_2\leqslant x/c_2^6\\ n_2\equiv n_1c_1^6c_2^{-6}(q)||\\  n_2 \equiv \alpha n_1c_1^6c_2^{-6}(q)}}r(n_2)
    \\ &+O\left(\sum_{b\leqslant \frac{x^{1/3}}{10c_1^2q^{2/3}}}\sum_{a\in[\frac{x^{1/2}}{c_1^3b^{3/2}}-q, \frac{x^{1/2}}{c_1^3b^{3/2}}]}\sum_{\substack{n_2\leqslant x/c_2^6\\ n_2\equiv n_1c_1^6c_2^{-6}(q)||\\  n_2 \equiv \alpha n_1c_1^6c_2^{-6}(q)}} |r(n_2)|\right),
\end{align*}

where the prime sum indicates that the number of $a$ is divisible by $q$ since we want to exploit circulating $n_2$ to traverse all elements of $(\Z/q\Z)^\times$. By Equation \eqref{eq:err-23-div} summing along $c_1,c_2$, the main term is now 
\begin{align*}
    \sum_{c_1,c_2\leqslant x^{1/6}}\sum_{b\leqslant \frac{x^{1/3}}{10c_1^2q^{2/3}}}\frac{x^{1/2}}{qc_1^3b^{3/2}}\sum_{n\leqslant x/c_2^6} r(n)\ll \dfrac{x^{\frac{1}{2}+\frac{2}{15}+\frac{1}{6}+\frac{1}{30}}}{q} \ll x^{1/6}q^{5/6}(x/q)^{-\eps}
\end{align*}
and the last inequality holds for all $q\gg x^{4/11+5\eps}$. By using the trivial bound \eqref{eq:trivialboundon_r(n)} on the error terms, we have that it contributes
\begin{align*}
    \sum_{c_1,c_2\leqslant x^{1/6}}\frac{x^{\frac{1}{3}}}{c_1^2q^{\frac{2}{3}}}\cdot q\frac{x^\frac{1}{2}}{c_2^3q}\ll \frac{x^{\frac{5}{6}}}{q^{\frac{2}{3}}} \ll x^{1/6}q^{5/6}(x/q)^{-\eps/10} 
\end{align*}
whenever $q\geqslant x^{4/9+\eps}$. Now, we bound the case when $b\gg \frac{x^{1/3}}{c_1^2q^{2/3}}$. This contributes to at most
\begin{align*}
    \sum_{b\geqslant \frac{x^{1/3}}{10c_1^2q^{2/3}}}&\sum_{a\leqslant \frac{x^{1/2}}{c_1^3b^{3/2}}}\sup_{\ell\in(\Z/q\Z)^\times} \sum_{\substack{n_2\leqslant x/c_2^6\\n_2\equiv \ell (q)}} |r(n_2)|
    \\ &\ll\sum_{b\geqslant \frac{x^{1/3}}{10c_1^2q^{2/3}}} \dfrac{x^{\frac{1}{2}}}{c_1^3b^{\frac{3}{2}}}\cdot \dfrac{x^{\frac{1}{2}}}{qc_2^3}\ll \dfrac{1}{c_1^3c_2^3}\dfrac{x}{q}\cdot \dfrac{c_1q^{\frac{1}{3}}}{x^{\frac{1}{6}}}=\dfrac{1}{c_1^2c_2^3}\cdot\frac{x^{\frac{5}{6}}}{q^{\frac{2}{3}}}.
\end{align*}
Summing along $c_1,c_2$ we see that this error term is acceptable when $q\geqslant x^{\frac{4}{9}+\eps}.$ We are done proving Proposition \ref{prop:err-term}.
\end{proof}
\subsection{Deduction of Theorem \ref{thm:sqfull-in-arith}}
Firstly, we determine the main term occurring from square-full integers that are relatively prime to $q$. We notice that
\begin{align*}
    \dfrac{1}{\varphi(q)} &\sum_{\substack{n\leqslant x\\ (n,q)=1}}1^\blacksquare(n) = \dfrac{1}{\varphi(q)}\left(\sum_{\substack{n\leqslant x}} 1^\blacksquare(n)-\sum_{\substack{n\leqslant x\\ q|n}} 1^\blacksquare(n)\right)\\
    &= \dfrac{\zeta(\frac{3}{2})}{\zeta(3)}x^{\frac{1}{2}}\frac{1}{q}\left(1-\frac{1}{q}\right)+\dfrac{\zeta(\frac{2}{3})}{\zeta(2)}x^{\frac{1}{3}}\frac{1}{q}\left(1-\dfrac{1}{q^{\frac{2}{3}}}\right)+O\left(\dfrac{x^{\frac{1}{6}}}{q}+\dfrac{x^{\frac{1}{2}}}{q^{\frac{5}{2}}}\right),
\end{align*}
where we have been using Equation \eqref{eq:counting-sqfull}. The last error term occurs from the second sum bounding $q$ multiple of square-full integers not exceeding $x/q^2$. The error term is acceptable when $q\gg x^{\frac{5}{29}+\eps}$.\\

This calculation with the consequence from Propositions \ref{prop:equiv-to-gcd-smallc} -- \ref{prop:err-term}, i.e., Equations \eqref{eq:eq-combine-4.1,4.2} and \eqref{eq:eq-4.5}, and the Cauchy-Schwarz inequality yield 
\begin{align*}
    \dfrac{1}{\varphi(q)} &\sum_{\ell\in (\Z/q\Z)^\times} \left| \frac{1}{2}\sum_{\substack{n\equiv \ell,\alpha\ell(q)\\ x<n\leqslant 2x}} 1^{\blacksquare}(n)-\left(\dfrac{\zeta(\frac{3}{2})}{\zeta(3)}\frac{1}{q}\left(1-\frac{1}{q}\right)\mathscr D_{\frac{1}{2}}(x;x)+\dfrac{\zeta(\frac{2}{3})}{\zeta(2)}\frac{1}{q}\left(1-\dfrac{1}{q^{\frac{2}{3}}}\right)\mathscr{D}_{\frac{1}{3}}(x;x)\right)
         \right|^2 \nonumber \\
         &\sim C(x/q)^{\frac{1}{6}},
\end{align*}
for all $x^{51/114+\eps}\leqslant q\leqslant x^{1-50\eps}$, which is our Theorem \ref{thm:sqfull-in-arith}.

\printbibliography[
heading=bibintoc,
title={References}
]

@misc{Kratzel-theorem,
    key="Kratzel",
    note={Kr\"{a}tzel, E.: Ein Teilerproblem (German), \textit{J. Reine Angew. Math.} \textbf{235} (1969) 150--174.}
}

@misc{error-term-in-square-full,
    key="Balasubramanian",
    note={R. Balasubramanian, K. Ramachandra, M. V. Subbarao, On the error function in the asymptotic formula for the counting function of $k$-full numbers, \textit{Acta Arith.} \textbf{50} (1988) 107--118.}
}

@misc{Bateman-Grosswald,
    key="Bateman",
    note={P. T. Bateman and E. Grosswald, On a theorem of Erd\H{o}s and Szekeres, \textit{Illinois J. Math.} \textbf{2} (1958) 88--98.}
}

@misc{square-full-Erdos-Szekeres,
    key="Erd\H{o}s",
    note={P. Erd\H{o}s and G. Szekeres, \"{U}ber die Anzahl der Abelschen Gruppen gegebener Ordnung und \"{u}ber ein verwandtes zahlentheoretisches Problem (German), \textit{Acta Sci. Math. (Szeged)} \textbf{7} (1934--1935) 95--102.}
}

@misc{Tenenbaum,
    key="Tenenbaum",
    note = {G. Tenenbaum, \textit{Introduction to Analytic and Probabilistic Number Theory, 3rd ed.}, AMS, 2015}
}

@misc{Titchmarsh,
    key="Titchmarsh",
    note = {E. C. Titchmarsh and D. R. Heath-Brown, \textit{The Theory of the Riemann Zeta-function}, OUP, 1986}
}

@misc{GoroMato,
    key="Gorodetsky",
    note={O. Gorodetsky, K. Matom\"{a}ki, M. Radziwi\l\l, and B. Rodgers, On the variance of squarefree integers in short intervals and arithmetic progressions, \textit{Geom. Func. Anal.} \textbf{31} (2021) 111--149.}
}

@misc{Chan-arith,
    key="Chan",
    note={T. H. Chan, Squarefull numbers in arithmetic progression II, \textit{J. Number Theory} \textbf{152} (2015) 90--104.}
}

@misc{chan-tsang,
    key="Chan",
    note={T. H. Chan, and K. M. Tsang, Squarefull numbers in arithmetic progression, \textit{Int. J. Number Theory} \textbf{9} (2013) 885--901.}
}

@misc{Chan1,
    key="Chan",
    note={T. H. Chan, Variance of squarefull numbers in short intervals, \textit{Illinois J. Math.} \textbf{67} (2023) 789--808.}
}

@misc{Chan2,
    key="Chan",
    note={T. H. Chan, Variance of squarefull numbers in short intervals II, \textit{arXiv preprint} arXiv:2311.13463}
}

@misc{Lester,
    key="Lester",
    note={S. Lester, On the variance of sums of divisor functions in short intervals, \textit{Proc. Amer. Math. Soc.} \textbf{144} (2016) 5015--5027}
}

@misc{Ivic-divisor-general,
    key="Ivic",
    note={A. Ivi\'{c}, The general divisor problem, \textit{J. Number Theory} \textbf{27} (1987) 73--91}
}

@misc{Ivic-divisor-short,
    key="Ivic",
    note={A. Ivi\'{c}, On the divisor function and the Riemann zeta-function in short intervals, \textit{Ramanujan J.} \textbf{19} (2009) 207--224}
}

@misc{Suryanarayana-generalized-sqfull,
    key="Suryanarayana",
    note={D. Suryanarayana, On the distribution of some generalized square-full integers, \textit{Pac. J. Math.} \textbf{72} (1977) 547--555}
}

@misc{Wang,
    key="Wang",
    note={D. Wang, On the distribution of square-full integers, \textit{Indian J. Pure Appl. Math.} \textbf{53} (2022) 627--634}
}

@misc{Trifonov,
    key="Trifonov",
    note={O. Trifonov, Lattice points close to a smooth curve and squarefull numbers in short intervals, \textit{J. London Math. Soc.} \textbf{65} (2002) 303--319}
}

@misc{Saffari-Vaughan,
    key="Saffari",
    note={B. Saffari and R. C. Vaughan. On the fractional parts of $x/n$ and related sequences. II, \textit{Ann. Inst. Fourier (Grenoble)} \textbf{27} (1977) 1--30}
}

@misc{Petrow-Young,
    key="Petrow",
    note={I. Petrow and M. P. Young. The fourth moment of Dirichlet L-functions along a coset and the Weyl bound, \textit{Duke Math. J.} \textbf{172} (2023) 1876--1960}
}

@misc{Montgomery,
    key="Montgomery",
    note = {H. L. Montgomery, \textit{Topics in multiplicative number theory}, Lecture Notes in Mathematics, Springer-Verlag, Berlin-New York, 1971}
}

@misc{Iwaniec-Kowalski,
    key="Iwaniec",
    note = {H. Iwaniec and E. Kowalski, \textit{Analytic number theory}, volume 53. American Mathematical Soc., 2004.}
}

\end{document}